\documentclass[11pt,preprint]{elsarticle}
\usepackage[a4paper,total={6.5in, 8in}]{geometry}

\makeatletter
\def\ps@pprintTitle{%
  \let\@oddhead\@empty
  \let\@evenhead\@empty
  \let\@oddfoot\@empty
  \let\@evenfoot\@oddfoot
}
\makeatother

\usepackage{pgfplots}
\pgfplotsset{compat=1.16} 

\usepackage{tikz}
\usepackage{amsmath,amssymb}
\usepackage{color}
\usepackage{pgfplots}
\usepackage{algorithm}
\usepackage{algpseudocode}
\usepackage{subfig}
\usepackage{graphicx}
\usetikzlibrary{plotmarks}
\usepackage{eufrak}
\usepackage{enumitem}

\usepackage{comment}

\usepackage{amsthm}
\theoremstyle{plain}
\newtheorem{theorem}{Theorem}[section]

\theoremstyle{definition}

\newtheorem{remark}[theorem]{Remark}

\theoremstyle{plain}

\newcommand{\beq}{\begin{equation}}
\newcommand{\eeq}{\end{equation}}
\newcommand{\bec}{\begin{center}}
\newcommand{\eec}{\end{center}}

\newcommand{\der}{\partial}

\newcommand{\Bb}{\mathbf{B}} 
\newcommand{\Cb}{\mathbf{C}}

\newcommand{\Gb}{\mathbf{G}}

\newcommand{\Jb}{\mathbf{J}}
\newcommand{\Kb}{\mathbf{K}}

\newcommand{\Ub}{\mathbf{U}}
\newcommand{\Vb}{\mathbf{V}}

\newcommand{\fb}{\mathbf{f}} 
\newcommand{\gb}{\mathbf{g}} 
 
\newcommand{\ib}{\mathbf{i}}

\newcommand{\nb}{\mathbf{n}}

\newcommand{\tb}{\mathbf{t}}
\newcommand{\ub}{\mathbf{u}}
\newcommand{\vb}{\mathbf{v}}

\newcommand{\Pib}{\boldsymbol{\Pi}}

\newcommand{\Acal}{\mathcal{A}}

\newcommand{\Kcal}{\mathcal{K}}
\newcommand{\Lcal}{\mathcal{L}} 
 
\newcommand{\Ncal}{\mathcal{N}}

\newcommand{\Pcal}{\mathcal{P}}

\newcommand{\Tcal}{\mathcal{T}}

\newcommand{\Vcal}{\mathcal{V}}

\newcommand{\Zcal}{\mathcal{Z}}

\newcommand{\iu}{\underline{i}}
\newcommand{\ju}{\underline{j}}

\newcommand{\abs}[1]{\left|#1\right|}
\newcommand{\ABS}[1]{\big|#1\big|}

\newcommand{\eqdef}{:=}
\newcommand{\defeq}{=:}
\newcommand{\dS}{dS}

\newcommand{\JUMP}[2]{\left[\!\left[#1\right]\!\right]_{#2}}
\newcommand{\restrict}[2]{{#1}_{|#2}}

\newcommand{\QTT}{QTT}

\newcommand{\StfMat}{\Kb}

\newcommand{\refsol}{{\textit{ref}{\hspace{0.5mm}}}}

\biboptions{sort&compress}

\begin{document}
\begin{frontmatter}

  \title{A Low-Rank QTT-based Finite Element Method for Elasticity
    Problems}

  \author[1]{Elena Benvenuti}
  \ead{elena.benvenuti@unife.it}
  
  \author[2]{Gianmarco Manzini\corref{cor2}}
  \ead{gmanzini@lanl.gov}
  \cortext[cor2]{Corresponding author.}
  
  \author[1]{Marco Nale}
  \ead{marco.nale@unife.it}
  
  \author[1]{Simone Pizzolato}
  \ead{simone.pizzolato@edu.unife.it}
  
  \address[1]{Department of Engineering, University of Ferrara, via
    Saragat 1, Ferrara, 44122, Italy}

  \address[2]{Theoretical Division, Los Alamos National Laboratory,
    Los Alamos, NM 87545, United States of America}
  
            
  
  \begin{abstract}
    We present an efficient and robust numerical algorithm for solving
    the two-dimensional linear elasticity problem that combines the
    Quantized Tensor Train format and a domain partitioning strategy.
    This approach makes it possible to solve the linear elasticity
    problem on a computational domain that is more general than a
    square.
    Our method substantially decreases memory usage and achieves a
    notable reduction in rank compared to established Finite Element
    implementations like the FEniCS platform.
    This performance gain, however, requires a fundamental rethinking
    of how core finite element operations are implemented, which
    includes changes to mesh discretization, node and degree of
    freedom ordering, stiffness matrix and internal nodal force
    assembly, and the execution of algebraic matrix-vector operations.
    In this work, we discuss all these aspects in detail and assess
    the method's performance in the numerical approximation of three
    representative test cases.
  \end{abstract}
  
  \begin{keyword}
    Finite element method         \sep
    tensor train format           \sep
    quantized tensor train format \sep
    elasticity    
  \end{keyword}
  
\end{frontmatter}




\section{Introduction}
\label{sec1:introduction}

The numerical approximation of partial differential equations (PDEs)
is pivotal in various scientific and engineering disciplines, for
example, in structural mechanics, especially when dealing with complex
geometries and singular solutions.
Traditional numerical methods, such as the Finite Element Method
(FEM), are widely recognized for their robustness and accuracy in
handling such problems.
However, the computational cost of these methods escalates rapidly
with increasing problem size and complexity, often rendering
high-resolution simulations impractical on conventional hardware.

Several strategies have been developed to address computational
resource limitations and efficiency.
Order reduction methods such as Ladèveze's LATIN method
\cite{ladeveze2010,ladeveze2012} and the Proper Generalized
Decomposition \cite{chinesta2013} have shown a significant performance
improvement by reducing the problem's dimensionality while maintaining
acceptable accuracy.
For instance, in computational homogenization, a computational
complexity of $\mathcal{O}\left(N^d\right)$ is usually required, with
$N$ being the number of degrees of freedom in each spatial dimension
(or a convenient upper bound of it) and $d$ being the number of
dimensions.
This exponential scaling with $d$ is called the \emph{curse of
dimensionality}~\cite{Bellman1961}, and is a major issue in the
numerical resolution of high-dimensional problems.


The development of data compression techniques based on tensor-network
formats, particularly the tensor train (TT)
format~\cite{oseledets2011,oseledets2010a}, has achieved a
breakthrough in high-dimensional data compression, efficient numerical
computations, and multi-dimensional array representation.
The tensor train format, which emerged from intensive research
\cite{oseledets2009b,oseledets2011,oseledets2009new,oseledets2009recursive},
offers an effective strategy for handling complex tensorial
structures.
This same mathematical framework, independently developed in quantum
physics
\cite{schollwock2011density,verstraete2006matrix,vidal2003efficient,white1993density},
is known as Matrix Product States (MPS) and serves as a fundamental
tool for analyzing quantum spin systems.

The distinguishing characteristics of the TT framework make it
exceptionally valuable for solving partial differential equations
(PDEs).
A key advantage lies in its storage efficiency: while traditional
tensor representations face exponential growth in memory requirements
with increasing dimensionality, the TT format achieves linear scaling
in the number of separated dimensions, with only quadratic dependence
on the ranks.
Therefore, the computational complexity scales like
$\mathcal{O}\left(dr^2N\right)$ instead of
$\mathcal{O}\left(N^d\right)$, where $r$ is a suitable upper bound on
the TT ranks, and it is significantly reduced whenever $r\ll N$, i.e.,
when a \emph{low-rank approximation} is feasible.
Additionally, it is worth emphasizing that a tensor train
decomposition algorithm exists that enables the computation of
``quasi-optimal'' approximations of any given tensors through a
systematic application of singular value decompositions on auxiliary
matrices, see, e.g.,
\cite[Alg.~1]{oseledets2011}.
The interplay between the accuracy of such approximation and the
computational costs can be controlled directly through a user-defined
tolerance factor $\epsilon$.

The practical utility of the TT format in designing efficient and
effective algorithms that approximate partial differential equations
is enhanced by its compatibility with fundamental linear algebraic
operations.
The tensor train framework incorporates efficient mechanisms for these
computations, complemented by robust rank reduction procedures, such as
\cite{verstraete2006matrix}
and
\cite[Alg.~2]{oseledets2011}, for managing the complexity of
intermediate results.
These features collectively establish the TT format as a powerful tool
for tackling high-dimensional computational challenges.
The mathematical foundation of this approach ensures both numerical
stability and computational efficiency, making it particularly
valuable for advanced scientific computing applications.

The \emph{Quantized Tensor Train (QTT)} format \cite{khoromskij2011}
extends this approach by introducing the concept of
``\emph{quantization of tensor dimensions}''.
Adopting the QTT format may lead to enhanced compression rates and
efficiency, which is particularly effective in solving PDEs.
Although this mathematical concept is relatively new in the field of
low-rank representations, it has rapidly demonstrated its remarkable
utility across numerical linear algebra and computational science
applications.
In fact,
numerous discrete operators naturally exhibit low-rank structure when
expressed in the QTT format, such as
\cite{oseledets2009a,oseledets2010b,Grasedyck2010,kazeev2012,kazeev2013b,khoromskij2018},
and this characteristic proves fundamental in developing efficient
tensor-based numerical approximations to PDE solutions.
Furthermore, the QTT-FEM methodology, which integrates low-rank QTT
decomposition with a fine-grained low-order finite element method,
presents a robust framework for obtaining accurate solutions for
multidimensional PDEs while maintaining computational efficiency.
In the last decade, extensive research
\cite{Dolgov2012,kazeev2014,kazeev2018b,bachmayr2020,khoromskij2011,kazeev2017,kazeev2022,kazeev2013a,markeeva2021}
has explored various partial differential models using this approach.
Notably, investigations into elliptic PDEs featuring singularities or
high-frequency oscillatory behaviors
\cite{kazeev2018,bachmayr2020,khoromskij2011,kazeev2017,kazeev2022}
have demonstrated exponential convergence rates relative to the total
parameter count, achieving results comparable to those observed in
hp-FEM applications for singular solutions.
The framework's effectiveness is further enhanced through advanced
preconditioning strategies.
The development of tensor-structured BPX preconditioners, initially
proposed in
\cite{bachmayr2020}
for uniformly elliptic problems and later extended to one-dimensional
singularly perturbed scenarios
\cite{khoromskij2018}
enables the practical implementation of low-rank QTT representations
on extremely refined grids.
This capability to handle grid resolutions approaching machine
precision typically eliminates the need for adaptive mesh refinement
procedures.

Recently published scientific literature focused on the application of
QTT-based numerical methods across diverse PDE categories.
Second-order elliptic PDEs have received significant attention
\cite{kazeev2015,kazeev2018,khoromskij2011,kazeev2017,kazeev2022,khoromskij2018,markeeva2021}
while specific investigations targeted
the one-dimensional Helmholtz equation
\cite{kazeev2017,fraschini2024a},
the chemical master equation
\cite{kazeev2014,kazeev2015b},
the molecular Schrödinger equation
\cite{kazeev2015},
and Fokker-Planck equation
\cite{Dolgov2012}.
The work we present in this paper has mainly been motivated by the
work of Reference~\cite{markeeva2021}, where the application of the
QTT format is explored within a FEM framework to solve the
two-dimensional Poisson equation on polygonal domains.
Our work extends the approach of References~\cite{markeeva2021} to the
linear elasticity model.
As in~\cite{markeeva2021}, we highlight the potential of the QTT
format to reduce memory consumption and to improve computational speed
compared to traditional sparse matrix representations, particularly
for fine meshes.
Moreover, using QTT in a domain-splitting setting makes it possible to
generalize the use of QTT format to computational domains with a more
general geometric shape, as, for example, the ``L-shape'' that we
consider in our numerical experiments.
Achieving optimal performance in the domain-splitting framework
demands a comprehensive redesign of how a finite element solver works,
with special care to the mesh organizations, a different reordering of
nodes and degrees of freedom, restructuring of stiffness matrix, and
internal nodal force assembly procedures, and, finally, a
reimplementation of basic algebraic operations.
A crucial point in keeping the possible rank growth under control is
the adoption of a special renumbering of the degrees of freedom that
follows the so-called \emph{Z-order}, see
~\cite{markeeva2020,morton1966}, and, accordingly, the introduction of
Z-kron operations~\cite{markeeva2021}, enabling the construction of
the stiffness matrix in Z-order directly in QTT format.
The main contribution of our work relies on a thorough discussion of
all the details and technicalities needed by the QTT-FEM and
experimentally proving its effectiveness in solving the linear
elasticity equation in variational form.
Our main result is that the combined strengths of FEM's adaptability
to complex geometries, and QTT's computational efficiency can
significantly enhance the performance of numerical simulations,
particularly in challenging scenarios involving singularities.


The paper is organized as follows.
Section \ref{sec2:math:model} introduces the mathematical formulation
of the linear elasticity problem, discussing both its strong and weak forms.
Section \ref{sec3:FEM} presents the finite element formulation in a domain
partitioning framework.
Section \ref{sec4:TT-QTT-formats} discusses a reformulation of the FEM 
by using the quantized tensor train format.
Section \ref{sec5:numerical} presents the result of our numerical experiments.
Section \ref{sec6:conclusions} offers final remarks and some hints about possible
future work.



\section{The linear elasticity model}
\label{sec2:math:model}

The linear elasticity model describes the small deformation behavior
of an elastic material under external forces or displacements.
We typically express the strong form of the linear elasticity problem
through the equilibrium equation with appropriate boundary conditions.

First, let $\Omega\subset\mathbb{R}^2$ be a two-dimensional domain
occupied by an elastic body, and $\partial\Omega$ its boundary, which
can be divided into two disjoint parts: $\partial\Omega_u$ where
displacements are prescribed, and $\partial\Omega_\sigma$ where
tractions (stresses) are prescribed.
Then, let $\ub=(u_{s})_{s=1,2}\in\mathbb{R}^2$ be the displacement field and
$\varepsilon(\ub)\in\mathbb{R}^{2\times 2}$, the strain, i.e., the
symmetric gradient
$\varepsilon(\ub) = \frac{1}{2} \left( \nabla\ub + (\nabla\ub)^T \right),$
where $\nabla\ub = \left(\partial u_s\slash{\partial x_t}\right)_{s,t=1,2}$.
Consider the stress tensor $\sigma(\ub)\in\mathbb{R}^{2\times2}$,
which is related to the strain tensor $\varepsilon(\ub)$ by the
constitutive (material) law, typically Hooke's law for linear elastic
materials, $\sigma(\ub) = C : \varepsilon(\ub),$ where $C$ is the
fourth-order elasticity tensor.
The equilibrium equation is:
\begin{align}
  \nabla\cdot\sigma + \fb =0 \quad \text{in }\Omega,
  \label{eq:strong:linear:elasticity}
\end{align}
where the vector-valued field $\fb$ represents the body forces.
To properly define the mathematical model, we need the boundary
conditions:
\begin{subequations}
  \label{eq:bnd}
  \begin{align}
    \ub            &= \ub_0\phantom{\tb} \, \text{on } \partial\Omega_u,    \label{eq:bnd:A}\\
    \sigma\cdot\nb &= \tb\phantom{\ub_0} \, \text{on } \partial\Omega_\sigma,\label{eq:bnd:B}
  \end{align}
\end{subequations}
where $\ub_0$ in ~\eqref{eq:bnd:A} is a prescribed displacement on the
Dirichlet boundary $\partial\Omega_u$; $\tb$ in~\eqref{eq:bnd:B} is
the prescribed traction on the Neumann boundary
$\partial\Omega_\sigma$; $\nb$ in~\eqref{eq:bnd:B} is the unit vector
orthogonal to $\partial\Omega_{\sigma}$ and pointing out of $\Omega$.
We assume that the boundaries $\partial\Omega_u$ and
$\partial\Omega_\sigma$ are nonoverlapping, in the sense that
$\abs{\partial\Omega_u\cap\partial\Omega_{\sigma}}=0$, where
$\abs{\,\cdot\,}$ is the one-dimensional Lebesgue measure of its
argument.

\medskip
As usual, we derive the weak form of
problem~\eqref{eq:strong:linear:elasticity}-\eqref{eq:bnd} by
multiplying~\eqref{eq:strong:linear:elasticity} by a vector-valued
test function $\vb$ that vanishes on $\partial\Omega_u$, and by
integrating by parts over the domain $\Omega$.
Using Hooke's law and noting that $C$ is a major symmetric tensor, the
variational form of the linear elasticity problem reads as:
\begin{align}
  &\emph{Find $\ub\in\Vb$ such that:}\nonumber\\
  & a(\ub,\vb) \eqdef
  \int_\Omega \nabla\ub : C : \nabla\vb\,d\Omega
  = \int_\Omega \fb\cdot\vb\,d\Omega + \int_{\partial\Omega_\sigma} \tb\cdot\vb\,\dS
  \defeq F(\vb)
  \quad \forall \vb\in\Vb_0,
  \label{eq:weak:general}
\end{align}
with the affine space
$\Vb=\big\{\vb\in\big[H^1(\Omega)\big]^2:\vb=\ub_0\text{~on~}\partial\Omega_u\big\}$,
and the linear subspace $\Vb_0$, which is obtained by setting
$\ub_0=0$ in $\Vb$.
The well-posedness of the variational
formulation~\eqref{eq:weak:general} follows from an application of the
Lax-Milgram theorem, which depends on the coercivity and continuity of
the bilinear form $a(\ub,\vb)$, and the continuity of the linear
functional $F(\vb)$.

This equation states that for all admissible test functions $\vb$, the
internal virtual work done by the stresses must equal the external
virtual work done by the body forces and tractions.
The elasticity tensor $C$ is a central component in the theory of
linear elasticity, representing the material's response to mechanical
stress.
Such a tensor plays a crucial role in linking the stress tensor
$\sigma$ to the strain tensor $\varepsilon$ through the Hooke's law,
which in tensor form is expressed as: $\sigma(\ub)=C\!:\!\varepsilon(\ub)$.
For isotropic materials, where properties are the same in all
directions, the tensor $C$ can be simplified significantly so that
it is characterized by only two independent constants, usually chosen
as the Young's modulus $E$ and the Poisson's ratio $\nu$.
The relationship, in this case, can be expressed as:
\begin{align*}
  C_{ijkl} = \lambda \delta_{ij} \delta_{kl}
  + \mu (\delta_{ik} \delta_{jl} + \delta_{il} \delta_{jk}),
\end{align*}
where $ \lambda $ and $\mu $ are the Lam\'e constants, related to $E$
and $\nu$ by
\begin{align*}
  \lambda = \frac{E \nu}{(1+\nu)(1-2\nu)}, \quad \mu = \frac{E}{2(1+\nu)}.
\end{align*}
The elasticity tensor is also related to the potential energy stored
in the material due to deformation.
Finally, the strain energy density $W$ for small deformations is given
by:
\begin{align*}
  W = \frac{1}{2} \varepsilon : C : \varepsilon.
\end{align*}
This relation expresses how energy is distributed and stored in the
material as a function of the strain, mediated by the properties
encapsulated in $C$.

\section{Finite element formulation}
\label{sec3:FEM}

To simplify the presentation, we consider the case where the boundary
of the domain associated with the stress variable,
$\partial\Omega_{\sigma}$, is the empty set so that the boundary of
the domain associated with the displacement variable,
$\partial\Omega_u$ coincides with the overall domain boundary,
$\partial\Omega$.
Additionally, we impose homogeneous Dirichlet boundary conditions by
setting to zero $\ub_0$, the prescribed displacement at the boundary
$\partial\Omega_u$.
Therefore, we solve the partial differential equation problem:
\begin{subequations}
  \begin{align}
    \nabla\cdot\sigma + \fb &= 0\phantom{0} \, \text{in }\Omega,        \label{eq:strong:A}\\
    \ub                     &= 0\phantom{0} \, \text{on }\partial\Omega \label{eq:strong:B}.
  \end{align}
\end{subequations}
The corresponding weak form reads as
\begin{align}
  &\emph{Find $\ub\in[H^1_0(\Omega)]^2$ such that:}\nonumber\\
  & \int_\Omega \sigma(\ub):\varepsilon(\vb)\,d\Omega
  = \int_\Omega \fb\cdot\vb\,d\Omega
  \quad \forall \vb\in[H^1_0(\Omega)]^2.
  \label{eq:weak:form}
\end{align}
Extending this formulation to the more general
formulation~\eqref{eq:weak:general} of the previous section is
straightforward, although requiring more technicalities.

To construct a finite element approximation of this problem that is
suitable to the QTT methodology, we consider the \emph{domain
partitioning technique}, following the idea originally proposed
in~\cite{markeeva2021} for the Poisson equation.
In practice, we split the open domain $\Omega$ into a set of $q$ open
subdomains $\Omega^{(m)}$ for $m=1,2,\ldots,q$, and
reformulate~\eqref{eq:weak:form} on each subdomain.
The regularity of the exact and approximate solutions implies that we
can ignore the interface integral terms that appear from a subdomain
integration by parts (a formal discussion about this point follows
below).
Then, we assemble the global stiffness matrix and the global load term
directly in the \QTT{} format from the local stiffness matrices and
load terms calculated on each subdomain, and we solve such a global
problem through the AMEn solver~\cite{dolgov2014}.
Hereafter, we present the domain partitioning approach for the finite
element discretization of~\eqref{eq:weak:form}; we will discuss the
application of the QTT methodology in the next section.

It is worth noting that \emph{our approach does not follow a domain
decomposition strategy}, where we should solve the partial
differential equation concurrently on each subdomain and reiterate
after exchanging information between the subdomains across their
common interfaces, see, cf.~\cite{dolean2015introduction}.
Nevertheless, a key point in this procedure is that the interface
nodes, which are shared by two adjacent subdomains, are replicated in
each subdomain in order to maintain $2^d$ nodes per side partition as
required by the QTT format.
Therefore, we will need a set of additional equations in the final
linear system to impose the numerical solution continuity at the
shared interface nodes.
This crucial point will be discussed in
Section~\ref{subsubsec:solution:concatenation}.

\subsection{Domain partitioning approach}
\label{subsec:DD:approach}

We partition the open domain $\Omega$ into $q$ open, quadrangular
subdomains $\Omega^{(m)}$ with boundary $\partial\Omega^{(m)}$, for
$m=1,2,\ldots,q$, such that
\begin{align*}
  \overline{\Omega} = \bigcup_{m=1}^{q} \overline{\Omega}^{(m)},
\end{align*}
where $\overline{\Omega}$ and $\overline{\Omega}^{(m)}$ denote the
closure in $\mathbb{R}^2$ of $\Omega$ and $\Omega^{(m)}$,
respectively.
We split the boundary of the $m$-th subdomain as
$\partial\Omega^{(m)} =
\Gamma^{(m,\text{ext})} \cup
\Gamma^{(m,\text{int})}$,
where
$\Gamma^{(m,\text{ext})}\eqdef\partial\Omega\cap\partial\Omega^{(m)}$
lies on the external boundary $\partial\Omega$, and
$\Gamma^{(m,\text{int})}\eqdef\partial\Omega^{(m)}\backslash\Gamma^{(m,\text{ext})}$
is the internal boundary, which has to be shared with at least another
subdomain.
From this definition, it trivially follows that the sub-boundaries
$\Gamma^{(m,\text{ext})}$ and $\Gamma^{(m,\text{int})}$ are
non-overlapping.
They can share mesh nodes, but they must also satisfy the condition
that
$\abs{\Gamma^{(m,\text{ext})}\cap\Gamma^{(m,\text{int})}}=0$, where
$\abs{\,\mathfrak{s}\,}$ is the one-dimensional Lebesgue measure in
$\mathbb{R}^2$ of set $\mathfrak{s}$ (the lenght of $\mathfrak{s}$).
To ease the notation, we will use the symbol $\Gamma^{(m)}$ instead of
$\Gamma^{(m,\text{int})}$ to denote the internal boundary of a single
domain $\Omega^{(m)}$.
Similarly, we will use the symbol
$\Gamma^{(m|p)}=\partial\Omega^{(m)}\cap\partial\Omega^{(p)}$ to
denote the interface specifically shared by the two subdomains
$\Omega^{(m)}$ and $\Omega^{(p)}$ (with $m\neq p$).

\medskip
After splitting the domain $\Omega$ into a set of $m$ disjoint
sub-domains $\Omega^{(m)}$, we consider a finite element approximation
of the partial differential equation
\begin{align}
  -\nabla\cdot(C:\nabla\ub^{(m)}) &= \fb^{(m)} \phantom{0} \,\, \text{in }\Omega^{(m)},                                 \label{eq:strong:Ab}\\[0.25em]
  \ub^{(m)}                       &= 0 \phantom{\fb^{(m)}} \,\, \text{on }\Gamma^{(m,\text{ext})}, \label{eq:strong:Bb}
\end{align}
where $\ub^{(m)}=\restrict{\ub}{\Omega^{(m)}}$,
$\fb^{(m)}=\restrict{\fb}{\Omega^{(m)}}$, and assuming that the flux
vector $\nb\cdot\sigma(\ub)$ is continuous at every internal
interface.
This continuity condition is consistent with the regularity assumption
that $\ub$ belongs to $[H^1_0(\Omega)]^2$.
Indeed, to derive the weak form
of~\eqref{eq:strong:Ab}-\eqref{eq:strong:Bb} on each subdomain
$\Omega^{(m)}$, we multiply both sides of \eqref{eq:strong:Ab} by the
vector-valued test function $\vb\in\big[H^1_0(\Omega)\big]^2$,
integrate by parts over $\Omega^{(m)}$ and add the resulting
expression over all the subdomains, so for $m=1,2,\ldots,q$.
We remove the external boundary integrals using \eqref{eq:strong:Bb}.
Then, let $\Omega^{m^+}$ and $\Omega^{m^-}$ be two subdomains sharing
a common interface
$\Gamma^{(m^+|m^-)}\supseteq\partial\Omega^{m^+}\cap\partial\Omega^{m^-}$.
We introduce the jump operator at $\Gamma^{(m^+|m^-)}$ that is defined
as
\begin{align*}
  \JUMP{C:\nabla\ub}{\Gamma^{(m^+|m^-)}} =
  \big[\nb\cdot(C:\nabla\ub)\big]^+ + \big[\nb\cdot(C:\nabla\ub)\big]^-,
\end{align*}
where the superscripts $\pm$ now indicate the side of the
$\Gamma$-interface on which the trace of the normal flux
$\nb\cdot(C\!:\!\nabla\ub)$ is taken.
After the integration by parts, we renumber all interfaces as
$\Gamma^{(s)}$, with $s$ ranging from $1$ to the total number of
interfaces.
Finally, we rearrange the interface summation term
$\Big[\sum_{m=1}^{q}\int_{\Gamma^{(m)}}\ldots\Big]$ as a summation
over all the subdomain interfaces, e.g.,
$\Big[\sum_{s}\int_{\Gamma^{(s)}}\ldots\Big]$, and we obtain:
\begin{align}
  &\sum_{m=1}^{q}\left(
  \int_{\Omega^{(m)}} \nabla\ub^{(m)}:C:\nabla\vb\,d\Omega -
  \int_{\Gamma^{(m)}} \big(\nb\cdot(C:\nabla\ub)\big)\vb\dS \right)
  \nonumber\\[0.5em]
  &\qquad=
  \sum_{m=1}^{q}\int_{\Omega^{(m)}} \nabla\ub^{(m)}:C:\nabla\vb\,d\Omega -
  \sum_{s}\int_{\Gamma^{(s)}}\JUMP{C:\nabla\ub}{\Gamma^{(s)}}\,\vb\,\dS
  \nonumber\\[0.5em]
  &\qquad=
  \int_{\Omega} \nabla\ub^{(m)}:C:\nabla\vb\,d\Omega.
  \label{eq:zero:jump}
\end{align}
The final equality in~\eqref{eq:zero:jump} is achieved since the
continuity of the flux implies that
$\JUMP{C:\nabla\ub}{\Gamma^{(s)}}=0$.

\subsection{Canonical and Z-ordering of the degrees of freedom}
\label{subsec:Z-ordering}
To define the finite element approximation in the subdomain
$\Omega^{(m)}$, it is essential to establish a mapping between the
nodal variables and a linear indexing system.
Specifically, let $\Omega_h^{(m)}$ be the mesh partitioning of the
subdomain $\Omega^{(m)}$, built as the tensor product of univariate
partitions with $n=2^d$ nodes per side.
We can identify every mesh node of $\Omega_h^{(m)}$ with the index
pair $(i,j)$, $i,j=0,\ldots,2^d-1$, and label each node with a single
index $\Lcal_{ij} = i + 2^d j$, ranging from $0$ to $4^{d}-1$.
We will refer to $\Lcal_{ij}$ as the \emph{canonical order} of the
nodes $\Omega^{(m)}$.

We will also consider a different enumeration called the
\emph{Z-order}.
This different numbering system identifies the $(i,j)$-th mesh node
with the number $\Zcal_{ij}$, which we uniquely determine by
interleaving the binary representation of indices $i$ and $j$,
see~\cite{morton1966}.
Hence, we first expand $i$ and $j$ on a binary digit format $[\ldots]_{2}$,
\begin{align*}
  i = \sum_{k=1}^d 2^{k-1} i_k \leftrightarrow i\equiv\big[i_1, i_2, \ldots i_d\big]_2
  \qquad\textrm{and}\qquad
  j = \sum_{k=1}^d 2^{k-1} j_k \leftrightarrow j\equiv\big[j_1, j_2, \ldots j_d\big]_2,
\end{align*}
where $i_d$ and $j_d$ are the most significant digits of $i$ and $j$,
respectively.
Then, we compute 
\begin{align*}
  \Zcal_{ij}
  &= i_1+2 j_1+4 i_2+8 j_2+\ldots+2^{2 d-2} i_d+2^{2 d-1} j_d
  = \sum_{k=1}^d 2^{2 k-2} i_k + \sum_{k=1}^d 2^{2 k-1} j_k \\ 
  &\leftrightarrow
  \Zcal_{ij} \equiv \big[ i_1, j_1, i_2, j_2, \ldots, i_d, j_d \big]_2.
\end{align*}

In the case of scalar fields, every mesh node is associated with a
nodal value, and, thus, with the number $\Lcal_{ij}$ or $\Zcal_{ij}$.
In the context of two-dimensional vector fields, we have two nodal
values associated with every mesh node, e.g.,
$\ub(x_i,y_j)=\big(u_x(x_i,y_j),u_y(x_i,y_j)\big)^T$ for the node with
coordinates $(x_i,y_j)$.
Consequently, we introduce a two-dimensional enumeration of the nodal
degrees of freedom for the canonical order,
\begin{align*}
  \Ncal_{ij}
  = \left[\begin{array}{c} \Ncal^{x}_{ij}  \\[0.5em]  \Ncal^{y}_{ij} \end{array}\right]
  = \left[\begin{array}{c}2\Lcal_{ij}+1   \\[0.5em] 2\Lcal_{ij}+2   \end{array}\right],
\end{align*}
and for the Z-order,
\begin{align*}
  \Pcal_{ij}
  = \left[\begin{array}{c}\Pcal^{x}_{ij} \\[0.5em] \Pcal^{y}_{ij}  \end{array}\right]
  = \left[\begin{array}{c}2\Zcal_{ij}+1  \\[0.5em] 2\Zcal_{ij}+2 \end{array}\right].
\end{align*}

\subsection{Finite element approximation on subdomain $\Omega^{(m)}$}
\label{subsec:FEM:subdomain}
For exposition's sake, we use hereafter the canonical order using
$\Lcal_{ij}$ and $\Ncal_{ij}$.
Switching to the Z-order formulation is almost straightforward and
implies using $\Zcal_{ij}$ and $\Pcal_{ij}$ instead of $\Lcal_{ij}$
and $\Ncal_{ij}$ in the following formula derivations.

On every subdomain $\Omega^{(m)}$, we consider the finite element
approximation of the displacement vector field $\ub^{(m)}$, which is
given by:
\begin{align*}
  \ub^{(m)}(x,y) \approx \sum_{i,j=1}^{2^d} \phi^{(m)}_{\Lcal_{ij}}(x,y) \Ub^{(m)}_{\Ncal_{ij}},
\end{align*}
where $\phi^{(m)}_{\Lcal_{ij}}(x,y)$ is the finite element shape
function associated with the local mesh node $\Lcal_{ij}$, and the
vector coefficients $\Ub^{(m)}_{\Ncal_{ij}}$ denotes the degrees of freedom associated with such a node.

We use the Lagrangian isoparametric elements; hence, the shape
function $\phi^{(m)}_{\Lcal_{ij}}$ is the piecewise bilinear function
that is one at node $\Lcal_{ij}$ and zero at all other nodes within
$\Omega_h^{(m)}$.
We refer to~\cite{bathe2006} for details on the
approximation properties of this discretization.
The displacement vector components in the $x$ and $y$ directions are expressed as:
\begin{align*}
  u^{(m)}_x(x,y) = \sum_{i,j=1}^{2^d}\phi^{(m)}_{\Lcal_{ij}}(x,y)\,U^{(m)}_{\Ncal^{x}_{ij}}
  \quad\textrm{and}\quad
  u^{(m)}_y(x,y) = \sum_{i,j=1}^{2^d}\phi^{(m)}_{\Lcal_{ij}}(x,y)\,U^{(m)}_{\Ncal^{y}_{ij}}.
\end{align*}
The array $\Ub^{(m)}$ holds the nodal unknowns and is organized as:
\begin{align*}
  \Ub^{(m)}
  = \Big(U^{(m)}_{\Ncal^{x}_{ij}}, U^{(m)}_{\Ncal^{y}_{ij}}\Big)^T.
\end{align*}
This leads to the local linear system:
\begin{align*}
  \StfMat^{(m)} \Ub^{(m)} + \big[\mbox{\text{interface terms}}\big] = \fb^{(m)},
\end{align*}
where
$\fb^{(m)}=\Big(f^{(m)}_{\Ncal^{x}_{ij}},f^{(m)}_{\Ncal^{y}_{ij}}\Big)^T$ is
the vector of nodal forces, defined by the subdomain integrals:
\begin{align*}
  f^{(m)}_{\Ncal^{x}_{ij}} = \int_{\Omega^{(m)}} f_x \phi^{(m)}_{\Lcal_{ij}} \,dx\,dy
  \quad\textrm{and}\quad
  f^{(m)}_{\Ncal^{y}_{ij}} = \int_{\Omega^{(m)}} f_y \phi^{(m)}_{\Lcal_{ij}} \,dx\,dy,
\end{align*}
and $\StfMat^{(m)}$ represents the stiffness matrix that is detailed
in the next subsections.
We left the interface terms unspecified because they add to zero when
we assemble the global linear system as proved in~\eqref{eq:zero:jump}.

\subsubsection{Solution concatenation}
\label{subsubsec:solution:concatenation}
We extend the procedure proposed by Markeeva et
al. \cite{markeeva2021} to impose a compatibility constraint between
the solutions of the system of equations on the subdomains to the
vector case.
Each subdomain is a quadrangle meshed with $2^d\times 2^d$ nodes.
Let $i,j \in [1\ldots 2^d]$ indicate the position of the generic
$(i,j)$-th node in the mesh.
In the two-dimensional elasticity vector case, we need to introduce
the multi-index $\Ncal_{ij}$ collecting the local index numbers in $x$
and $y$, i.e., $\Ncal^x_{ij}$ and $\Ncal^y_{ij}$.
Then, by using this notation and omitting the interface terms, we can
rewrite the local linear system of equations at node $\Ncal_{ij}$ as:
\begin{equation}
  \label{eq:equil_subdomain}
  \sum_{k,l =1} ^{2^d} \Kb^{(m)}_{\Ncal_{ij},\Ncal_{kl}} \Ub^{(m)}_{\Ncal_{kl}} = \fb^{(m)}_{\Ncal_{ij}}.
\end{equation}
When a node is located at the crossroad among $N_s$ subdomains, the
equilibrium equation \eqref{eq:equil_subdomain} can be rewritten
taking into account the additive contribution of these subdomains as
follows:
\begin{equation}
  \label{eq:equil_subdomain_junction}
  \sum_{s=1}^{N_s}\sum_{k,l=1}^{2^d} \Kb^{(s)}_{\Ncal_{ij},\Ncal_{kl}} \Ub^{(s)}_{\Ncal_{kl}}
  = \sum_{s=1}^{N_s}\fb^{(s)}_{\Ncal_{ij}}.
\end{equation}

For example, let subdomains $m_1$ and $m_2$ share a side, and assume
that the node $(i_1,\,j_1)$ from $m_1$ coincide with the node
$(i_2,\,j_2)$ from $m_2$.
The compatibility between the displacement at the junction node is
restored by imposing the following conditions:
\begin{subequations}
  \begin{align}\label{eq:compatibility}
    &\Ub^{(m_1)}_{\Ncal^x_{i_1 j_1}} =  \Ub^{(m_2)}_{\Ncal^x_{i_2 j_2}}\,\\
    &\Ub^{(m_1)}_{\Ncal^y_{i_1 j_1}} =  \Ub^{(m_2)}_{\Ncal^y_{i_2 j_2}}.
  \end{align}
\end{subequations}
Then, conditions \eqref{eq:compatibility} are satisfied via a
Lagrangian multiplier approach as follows:
\begin{align*}
  &\lambda
  \big(\Ub^{(m_1)}_{\Ncal^x_{i_1 j_1}} - \Ub^{(m_2)}_{\Ncal^x_{i_2 j_2}}\big) +
  \sum_{k,l=1}^{2^d} \Kb^{(m_1)}_{\Ncal_{i_1 j_1},\Ncal_{kl}} \Ub^{(m_1)}_{\Ncal_{kl}} +
  \sum_{k,l=1}^{2^d} \Kb^{(m_2)}_{\Ncal_{i_2 j_2},\Ncal_{kl}} \Ub^{(m_2)}_{\Ncal_{kl}} =
  \fb^{(m_1)}_{\Ncal_{i_1 j_1}} + \fb^{(m_2)}_{\Ncal_{i_2 j_2}},\\
  &\lambda
  \big(\Ub^{(m_1)}_{\Ncal^y_{i_1 j_1}} - \Ub^{(m_2)}_{\Ncal^y_{i_2 j_2}}\big) +
  \sum_{k,l=1}^{2^d} \Kb^{(m_1)}_{\Ncal_{i_1 j_1},\Ncal_{kl}} \Ub^{(m_1)}_{\Ncal_{kl}} +
  \sum_{k,l=1}^{2^d} \Kb^{(m_2)}_{\Ncal_{i_2 j_2},\Ncal_{kl}} =
  \fb^{(m_1)}_{\Ncal_{i_1 j_1}} + \fb^{(m_2)}_{\Ncal_{i_2 j_2}}\Ub^{(m_2)}_{\Ncal_{kl}},\\
  &\lambda
  \big(\Ub^{(m_2)}_{\Ncal^x_{i_2 j_2}} - \Ub^{(m_1)}_{\Ncal^x_{i_1 j_1}}\big) +
  \sum_{k,l=1}^{2^d} \Kb^{(m_1)}_{\Ncal_{i_1 j_1},\Ncal_{kl}} \Ub^{(m_1)}_{\Ncal_{kl}} +
  \sum_{k,l=1}^{2^d} \Kb^{(m_2)}_{\Ncal_{i_2 j_2},\Ncal_{kl}} \Ub^{(m_2)}_{\Ncal_{kl}} =
  \fb^{(m_1)}_{\Ncal_{i_1 j_1}} + \fb^{(m_2)}_{\Ncal_{i_2 j_2}},\\
  &\lambda
  \big(\Ub^{(m_2)}_{\Ncal^y_{i_2 j_2}} - \Ub^{(m_1)}_{\Ncal^y_{i_1 j_1}}\big) +
  \sum_{k,l=1}^{2^d} \Kb^{(m_1)}_{\Ncal_{i_1 j_1},\Ncal_{kl}} \Ub^{(m_1)}_{\Ncal_{kl}} +
  \sum_{k,l=1}^{2^d} \Kb^{(m_2)}_{\Ncal_{i_2 j_2},\Ncal_{kl}} \Ub^{(m_2)}_{\Ncal_{kl}} =
  \bf^{(m_1)}_{\Ncal_{i_1 j_1}} + \fb^{(m_2)}_{\Ncal_{i_2 j_2}},
\end{align*}
where the Lagrangian multiplier $\lambda$ is a positive real number.

\subsubsection{Construction of the stiffness matrix}
\label{subsubsec:stiffness:matrix}

We obtain the stiffness matrix
$\StfMat^{(m)}=\Big(\StfMat_{\Ncal_{ij}\Ncal_{kl}}^{(m)}\Big)_{\Ncal_{ij},\Ncal_{kl}=1,2,\ldots,4^{d}}$
from the weak formulation by restricting the integration from $\Omega$
in the bilinear form $a(\ub,\vb)$ to the subdomain $\Omega^{(m)}$.
The components associated with the connected nodes $\Lcal_{ij}$ and
$\Lcal_{kl}$ is:
\begin{align}
  \label{eq-stiffness}
  \StfMat_{\Ncal_{ij}\Ncal_{kl}}^{(m)}
  = \int_{\Omega^{(m)}} \Bb_{\mathcal{L}_{ij}}^T \Cb \Bb_{\mathcal{L}_{hk}} \, dx \, dy,
\end{align}
where $\Cb$ is the constant, constitutive, symmetric matrix, which for
the plane isotropic case is:
\begin{align}
  \Cb =
  \begin{bmatrix}
    C_{11} & C_{12} & 0 \\
    C_{12} & C_{22} & 0 \\
    0 & 0 & C_{33}
  \end{bmatrix},
  \label{eq:Cmat}
\end{align}
and $\Bb_{\Lcal_{ij}}$ is the compatibility operator at the local node
$(x_i, y_j)$ that can be derived from the strain-displacement
relationship
$\varepsilon_{\mathcal{L}_{ij}}=\Bb_{\mathcal{L}_{ij}}\Ub_{\Ncal_{ij}}$.
Recalling that $\Ub_{\Ncal_{ij}} = (U_{\Ncal^{x}_{ij}},
U_{\Ncal^{y}_{ij}})^T$, it holds that
\begin{align}
  \label{eq-B}
  \varepsilon(\ub^{(m)}) =
  \begin{bmatrix}
    \epsilon_{xx}(\ub^{(m)}) \\[0.5em]
    \epsilon_{yy}(\ub^{(m)}) \\[0.5em]
    2\epsilon_{xy}(\ub^{(m)})
  \end{bmatrix} =
  \begin{bmatrix}
    \displaystyle\frac{\partial u_x^{(m)}}{\partial x} \\[1em]
    \displaystyle\frac{\partial u_y^{(m)}}{\partial y} \\[1em]
    \displaystyle\frac{\partial u_x^{(m)}}{\partial y} + \frac{\partial u_y^{(m)}}{\partial x}
  \end{bmatrix}
  \quad\text{and}\quad
  \Bb_{\Lcal_{ij}} =
  \begin{bmatrix}
    \displaystyle\frac{\partial\phi^{(m)}_{\Lcal_{ij}}}{\partial x} & 0 \\
    0                                                           & \displaystyle\frac{\partial\phi^{(m)}_{\Lcal_{ij}}}{\partial y} \\[1em]
    \displaystyle\frac{\partial\phi^{(m)}_{\Lcal_{ij}}}{\partial y} & \displaystyle\frac{\partial\phi^{(m)}_{\Lcal_{ij}}}{\partial x}
  \end{bmatrix},
\end{align}
where we adopted the Voigt notation (as usual in finite element
discretization of elliptic problems \cite{bathe2006}).
Replacing \eqref{eq-B} and \eqref{eq:Cmat} in \eqref{eq-stiffness}, we
recast the latter as:
\begin{align}
  \StfMat_{\Ncal_{ij}\Ncal_{kl}}^{(m)} =
  \int_{\Omega^{(m)}}
  \begin{bmatrix}
    \displaystyle\frac{\partial\phi^{(m)}_{\Lcal_{ij}}}{\partial x} & 0                                                           & \displaystyle\frac{\partial\phi^{(m)}_{\Lcal_{ij}}}{\partial y} \\[1em]
    0                                                           & \displaystyle\frac{\partial\phi^{(m)}_{\Lcal_{ij}}}{\partial y} & \displaystyle\frac{\partial\phi^{(m)}_{\Lcal_{ij}}}{\partial x}
  \end{bmatrix}
  \,
  \begin{bmatrix}
    C_{11} & C_{12} & 0 \\
    C_{12} & C_{22} & 0 \\
    0     & 0     & C_{33}
  \end{bmatrix}
  \,
  \begin{bmatrix}
    \displaystyle\frac{\partial\phi^{(m)}_{\Lcal_{hk}}}{\partial x} & 0\\[1em]
    0                                                           & \displaystyle\frac{\partial\phi^{(m)}_{\Lcal_{hk}}}{\partial y} \\[1em]
    \displaystyle\frac{\partial\phi^{(m)}_{\Lcal_{hk}}}{\partial y} & \displaystyle\frac{\partial\phi^{(m)}_{\Lcal_{hk}}}{\partial x}
  \end{bmatrix}
  \,dx\,dy,
\end{align}
from which the following subdomain matrix descends:
\begin{align}
  \label{eq:stifblock}
  \StfMat_{\Ncal_{ij}\Ncal_{hk}}^{(m)} =
  \int_{\Omega}
  \begin{bmatrix}
    A_{\Ncal^{x}_{ij} \Ncal^{x}_{hk}}^{(m)} & A_{\Ncal^{x}_{ij} \Ncal^{y}_{hk}}^{(m)} \\[0.5em]
    A_{\Ncal^{y}_{ij} \Ncal^{x}_{hk}}^{(m)} & A_{\Ncal^{y}_{ij} \Ncal^{y}_{hk}}^{(m)}
  \end{bmatrix}
  \,dx\,dy.
\end{align}


\subsection{Subdomain numerical integration}
\label{subsec:subdomain:numerical:integration}
We remap each element of the subdomain mesh $\Omega^{(m)}_h$ onto the
adimensional reference element $\mathbb{K}=[-1,1]\times[-1,1]$ by
introducing a bijective mapping from the global coordinate system
$(x,y)$ to the dimensionless coordinate system $(\xi,\eta)$ defined on
$\mathbb{K}$.
Since we use the Lagrangian isoparametric elements, we consider the
four shape functions associated with the corners of $\mathbb{K}$,
which are
\begin{align*}
  \Phi_{-1,-1}(\xi,\eta)&=\dfrac{(1-\xi)(1-\eta)}{4},\quad
  \Phi_{ 1,-1}(\xi,\eta) =\dfrac{(1+\xi)(1-\eta)}{4},\\
  \Phi_{ 1, 1}(\xi,\eta)&=\dfrac{(1+\xi)(1+\eta)}{4},\quad
  \Phi_{-1, 1}(\xi,\eta) =\dfrac{(1-\xi)(1+\eta)}{4}.
\end{align*}
Let $(x_1,y_1)$, $(x_2,y_2)$, $(x_3,y_3)$, $(x_4,y_4)$ be the
coordinates of the four vertices of a given element of
$\Omega^{(m)}_h$, which are respectively mapped onto the four vertices
of $\mathbb{K}$ with coordinates $(-1,-1)$, $(1,-1)$, $(1,1)$,
$(-1,1)$.
To ease the notation, in this subsection we prefer dropping the
superindex $(m)$, so that we denote the coordinate mappings as
$(x_{\ell},y_{\ell})$ instead of $(x^{(m)}_{\ell},y^{(m)}_{\ell})$ for
$\ell=1,2,3,4$.
We similarly remove the superindex $(,)$ in all quantities related to
such mappings, e.g., their partial derivatives and Jacobian matrices.
The change of coordinates is given by
\begin{align*}
  \begin{bmatrix} x(\xi, \eta) \\ y(\xi,\eta) \end{bmatrix} =
  \begin{bmatrix}
    x_1 & x_2 & x_3 & x_4 \\
    y_1 & y_2 & y_3 & y_4 
  \end{bmatrix}\,
  \left[
    \begin{array}{l}
      \Phi_{-1,-1}(\xi,\eta) \\[0.5em]
      \Phi_{ 1,-1}(\xi,\eta) \\[0.5em]
      \Phi_{ 1, 1}(\xi,\eta) \\[0.5em]
      \Phi_{-1, 1}(\xi,\eta)    
    \end{array}
    \right],
\end{align*}
so that $\big(x(-1,1),y(-1,1)\big)=(x_1,y_1)$, etc.
The Jacobian matrix $\Jb(\xi,\eta)$ of the transformation from
$(\xi,\eta)$ to $(x,y)$ is non-singular.
Therefore, the inverse Jacobian matrix $\Jb^{-1}(x,y)$ associated with
the inverse mapping $\big(\xi(x,y),\eta(x,y)\big)$ from the reference
element $\mathbb{K}$ to a given mesh element of $\Omega^{(m)}_h$ is
also defined (and non-singular).
The Jacobian matrix and its inverse take the form:
\begin{align*}
  \Jb(\xi,\eta) =
  \begin{bmatrix}
    \dfrac{\der x}{\der\xi } & \dfrac{\der y}{\der\xi} \\[1em]
    \dfrac{\der x}{\der\eta} & \dfrac{\der y}{\der\eta}
  \end{bmatrix},
  \qquad
  \Jb^{-1}(\xi,\eta) =
  \dfrac{1}{|\Jb(\xi,\eta)|}
  \begin{bmatrix}
    \dfrac{\der\xi}{\der x} & \dfrac{\der\eta}{\der x} \\[1em]
    \dfrac{\der\xi}{\der y} & \dfrac{\der\eta}{\der y}
  \end{bmatrix},
\end{align*}
where $|\Jb(\xi,\eta)|$ is the determinant of $\Jb(\xi,\eta)$.

\medskip

Let $c\in\big\{ (-1,-1), (-1,1), (1,1), (-1,2)\big\}$ be a corner of
the reference element and $\Phi_{c}$ be the corresponding shape
function.
Then, consider a node $\Lcal_{ij}$ of the subdomain mesh
$\Omega^{(m)}_h$ and let $c$ be the corresponding node of $\mathbb{K}$
determined by the above mapping from a mesh element having the node
$\Lcal_{ij}$ as one of its nodes to the adimensional coordinate system
$(\xi,\eta)$.
Let $\Bb_{c}$ be the compatibility matrix $\Bb_{\Lcal_{ij}}$ at the
node $\Lcal_{ij}$ that is mapped to $c$ in $\mathbb{K}$ expressed in
the reference coordinate system $(\xi,\eta)$.
On applying the coordinate transformation, we find that
\begin{align*}
  \Bb_{c} &=
  \begin{bmatrix}
    \dfrac{\der\Phi_{c}}{\der\xi }\dfrac{\der\xi }{\der x} +
    \dfrac{\der\Phi_{c}}{\der\eta}\dfrac{\der\eta}{\der x} & 0 \\[1.5em]
    0 &
    \dfrac{\der\Phi_{c}}{\der\xi }\dfrac{\der\xi}{\der y} +
    \dfrac{\der\Phi_{c}}{\der\eta}\dfrac{\der\eta}{\der y} \\[1.5em]
    \dfrac{\der\Phi_{c}}{\der\xi }\dfrac{\der\xi }{\der y} +
    \dfrac{\der\Phi_{c}}{\der\eta}\dfrac{\der\eta}{\der y} &
    \dfrac{\der\Phi_{c}}{\der\xi }\dfrac{\der\xi }{\der x} +
    \dfrac{\der\Phi_{c}}{\der\eta}\dfrac{\der\eta}{\der x}
  \end{bmatrix}
  =
  \bigg(\dfrac{\der\Phi_{c}}{\der\xi},\dfrac{\der\Phi_{c}}{\der\eta}\bigg)
  \begin{bmatrix}
    \left(\dfrac{\der\xi }{\der x}, \dfrac{\der\eta}{\der x}\right)^T & 0 \\[1.5em]
    0 &
    \left(\dfrac{\der\xi}{\der y}, \dfrac{\der\eta}{\der y}\right)^T \\[1.5em]
    \left(\dfrac{\der\xi }{\der y}, \dfrac{\der\eta}{\der y}\right)^T &
    \left(\dfrac{\der\xi }{\der x}, \dfrac{\der\eta}{\der x}\right)^T
  \end{bmatrix}\\[1em]
  &=
  \nabla_{\xi,\eta}\Phi_c(\xi,\eta)\,\widehat{\Jb}(x(\xi,\eta),y(\xi,\eta)),
\end{align*}
where we have set
\begin{equation}\label{eq:Jhat}
  \widehat{\Jb}(x(\xi,\eta),y(\xi,\eta)) := \begin{bmatrix}
    \left(\dfrac{\der\xi }{\der x}, \dfrac{\der\eta}{\der x}\right)^T & 0 \\[1.5em]
    0 &
    \left(\dfrac{\der\xi}{\der y}, \dfrac{\der\eta}{\der y}\right)^T \\[1.5em]
    \left(\dfrac{\der\xi }{\der y}, \dfrac{\der\eta}{\der y}\right)^T &
    \left(\dfrac{\der\xi }{\der x}, \dfrac{\der\eta}{\der x}\right)^T
  \end{bmatrix}.
  \end{equation}  

The last equality above implicitly defines the matrix $\widehat{\Jb}$
that depends on the entries of $\Jb^{-1}$, i.e., on the derivatives of
$\xi(x,y)$ and $\eta(x,y)$ with respect to $x$ and $y$.
To ease the notation, let use $\widehat{\Jb}(\xi,\eta)$ to denote
$\widehat{\Jb}(x(\xi,\eta),y(\xi,\eta))$.
The contribution to the stiffness matrix from the pair of shape
functions $\Phi_{c_1}$ and $\Phi_{c_2}$ for $c_1,c_2\in\big\{(-1,-1),
(-1,1), (1,1), (1,-1)\big\}$ is
\begin{align}
  \Kb_{c_1,c_2}
  & = \int_{\mathbb{K}} \Bb^T_{c_1}(\xi,\eta)\Cb\Bb_{c_2}(\xi,\eta)\left|\Jb_{\xi \eta}\right|\,d\xi\,d\eta\nonumber\\[1.5em]
  & =
  \int_{\mathbb{K}} 
  \bigg[
    \widehat{\Jb}^T(\xi,\eta)\,
    \Big(\nabla_{\xi,\eta}\Phi_{c_1}(\xi,\eta)\Big)^T 
    \bigg]\,\,
  \Cb\,\,  
  \bigg[
    \nabla_{\xi,\eta}\Phi_{c_2}(\xi,\eta)\,
    \widehat{\Jb}(\xi,\eta)
    \bigg]
  \left|\Jb_{\xi \eta}\right| d \xi d \eta.
  \label{eq:Kc1c2:def}
\end{align}
The associated elastic energy is written as:
\begin{align*}
  E_{c_{1},c_{2}} = \dfrac{1}{2}\Ub_{c_1}^T \Kb_{c_1,c_2} \Ub_{c_2},
\end{align*}
where we denote $\Ub_{c_{1}}=\big(U_{c_{1,x}}, U_{c_{1,y}}\big)^T$ and
$\Ub_{c_{2}}=\big(U_{c_{2,x}}, U_{c_{2,y}}\big)^T$, according to
``$c$''-notation introduced above.
To evaluate the stiffness matrix integrals, we consider a Gaussian
quadrature rule that is obtained by the tensor product of the
one-dimensional quadrature rules with same weights and node
distributions $\{(W_{g_{i}},\xi_{g_{i}})\}$ and
$\{(W_{g_{j}},\eta_{g_{j}})\}$ along the directions $x$ and $y$,
respectively.
The numerical integration reads as:
\begin{align*}
  \Kb_{c_1,c_2} \approx
  \sum_{g_i=1}^{N_{G_{i}}}\sum_{g_j=1}^{N_{G_{j}}}
  W_{g_i}\,W_{g_j}
  \Bb^T_{c_{1}}(\xi_{g_i},\eta_{g_j})\Cb\Bb_{c_{2}}(\xi_{g_i},\eta_{g_j})
  \abs{\Jb(\xi_{g_i},\eta_{g_j}\eta)}.
\end{align*}
A further simplification is given by evaluating the integrand at the
element center.

\begin{figure}
    \centering
    \includegraphics[width=0.7\linewidth]{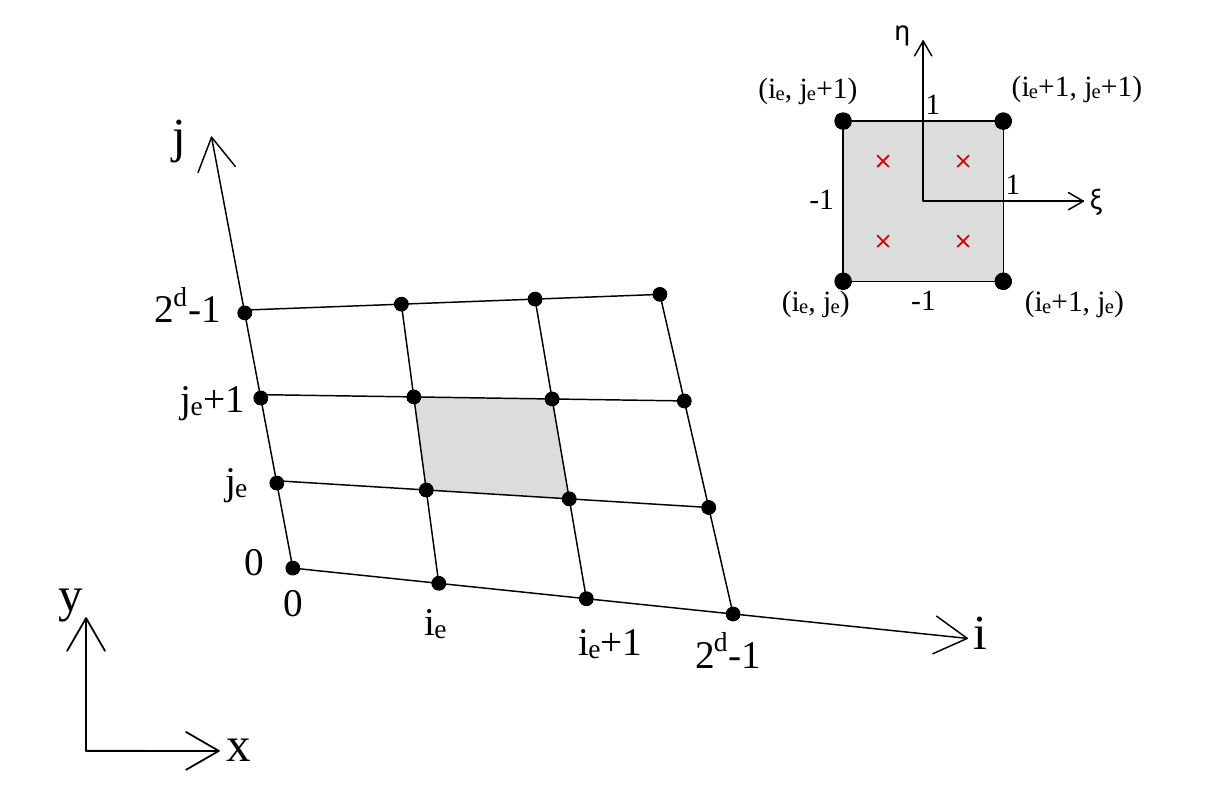}
    \caption{Mesh partition for the Jacobian calculation}
    \label{fig:jacobian}
\end{figure}

\medskip
We can also approximate the partial derivatives of the Jacobian matrix
in each mesh element of $\Omega^{(m)}_h$ as follows.
Consider the $(i,j)$-th element whose bottom-left corner is the mesh
node identified by the index pair $(i,j)$, with
$i,j=0,1,\ldots,2^d-1$.
We expand the entries of the Jacobian matrix in such an element,
denoted by $\Jb^{(i, j)}(\xi,\eta)$, as follows:
\begin{align}
  \Jb^{(i, j)}(\xi,\eta)
  = \Jb^{(0,0)}(\xi,\eta)
  + i\big( \Jb^{(1,0)}(\xi,\eta) - \Jb^{(0,0)}(\xi,\eta) \big)
  + j\big( \Jb^{(0,1)}(\xi,\eta) - \Jb^{(0,0)}(\xi,\eta) \big),
  \label{eq:Jij:expansion}
\end{align}
where $\Jb^{(0,0)}$, $\Jb^{(1,0)}$, and $\Jb^{(0,1)}$ are the Jacobian
matrices associated with the near elements that are identified by the
index pairs $(0,0),(1,0),(0,1)$, see Fig.~\ref{fig:jacobian}.
Using the Gauss nodes $(\xi_{g_i},\eta_{g_j})$, we rewrite the formula
above as:
\begin{multline*}
  \Jb^{(i, j)}(\xi_{g_i},\eta_{g_j})
  = \Jb^{(0,0)}(\xi_{g_i},\eta_{g_j})\\
  + i\big( \Jb^{(1,0)}(\xi_{g_i},\eta_{g_j}) - \Jb^{(0,0)}(\xi_{g_i},\eta_{g_j}) \big)
  + j\big( \Jb^{(0,1)}(\xi_{g_i},\eta_{g_j}) - \Jb^{(0,0)}(\xi_{g_i},\eta_{g_j}) \big).
\end{multline*}
Analogously, its determinant $\abs{\Jb^{(i,j)}(\xi,\eta)}$ is
calculated using the expression:
\begin{align}
  \ABS{\Jb^{(i,j)}(\xi,\eta)} =
  \ABS{\Jb^{(0,0)}(\xi,\eta)}+
  i\ABS{\Jb^{(1,0)}(\xi,\eta) - \Jb^{(0,0)}(\xi,\eta)} +
  j\ABS{\Jb^{(0,1)}(\xi,\eta) - \Jb^{(0,0)}(\xi,\eta)}.
\end{align}
Once $\Jb^{(i,j)}$ is computed, the operator
$\widehat{\Jb}^{(i,j)}(\xi,\eta)$ is straightforwardly obtained from
Eq.~\eqref{eq:Jhat}.


\subsection{Assembly at the subdomain level}
\label{subsec:assembly:subdomain:level}
All operations for computing $\Kb_{c_1,c_2}$ for a given pair
$(c_1,c_2)$ are vectorizable; so, we can carry out the calculation for
all the mesh elements of a particular domain, which are addressed by
$\Lcal_{ij}$, simultaneously.
Let $\Kcal_{(c_1,c_2)}$ be the array collecting all the elemental
integral values and $\Vb_{c_1}$ and $\Vb_{c_2}$ the shift matrices
moving such values to their final position into the stiffness
matrix $\StfMat^{(m)}_{(c_1c_2)}$ associated with the subdomain
$\Omega^{(m)}$.
We formally write that 
\begin{align}
  \StfMat^{(m)}_{(c_1,c_2)}
  = \Vb_{c_1}^T\text{diag}\big(\Kcal_{(c_1,c_2)}\big)\Vb_{c_2},
  \label{eq:Nodal:Stiffness:Matrix}
\end{align}
and we obtain the stiffness matrix $\StfMat^{(m)}$ of $\Omega^{(m)}$ by
accumulating the contributions from all the sixteen possible pairs
$(c_1,c_2)$, so that
\begin{align}
  \StfMat^{(m)} = \sum_{(c_1,c_2)} \StfMat^{(m)}_{(c_1,c_2)}.
  \label{eq:Subdomain:Stiffness:Matrix}
\end{align}
We can build matrix $\StfMat^{(m)}$ in canonical or Z-order by a
suitable choice of the shifting matrices $\Vb_{c_1}$ and $\Vb_{c_2}$.

Analogously, for mesh element of $\Omega^{(m)}_{h}$ and any one of the
sixteen possible combinations of $c_1$ and $c_2$, we compute the
right-hand side load vector $\Gb_{c_1,c_2}$ as
\begin{align*}
  \Gb_{c_1,c_2} =   \sum_{g_i=1}^{N_{G_{i}}}\sum_{g_j=1}^{N_{G_{j}}}
  W_{g_i}\,W_{g_j}
  \Phi_{c_1}(\xi_{g_i},\eta_{g_j})\Phi_{c_2}(\xi_{g_i},\eta_{g_j})
  \abs{\Jb(\xi_{g_i},\eta_{g_j})}
\end{align*}
where $\Jb_0$ is the Jacobian evaluated at the center of that mesh
element.
Then, we move the elemental values to their global position through
the same shift matrices $\Vb_{c_1}$ and $\Vb_{c_2}$, so that
\begin{align}
  \fb^{(m)}_{c_1,c_2} =
  \Vb_{c_1}^T\text{diag}(\Gb_{c_1,c_2})\Vb_{c_2}\overline{\fb},
  \label{eq:load:vector:intensity}
\end{align}
where $\overline{\fb}$ denotes the force vector intensity.
Finally, we consider the summation
\begin{align}
  \fb^{(m)} = \sum_{(c_1,c_2)} \fb^{(m)}_{c_1,c_2},
  \label{eq:subdomain:load:vector:intensity}
\end{align}
which accumulates all contributions~\eqref{eq:load:vector:intensity}
on the right-hand side (RHS) load vector $\fb^{(m)}$.

\subsection{Assembly at the global domain level}
\label{subsec:assembly:global:domain:level}
We split the stiffness matrix $\StfMat^{(m)}$ and the RHS load vector
$\fb^{(m)}$ of the subdomain $\Omega^{(m)}$ as follows:
\begin{align*}
  \StfMat^{(m)} =
  \begin{bmatrix}
    \Kb^{(m)}_{xx} & \Kb^{(m)}_{xy} \\
    \Kb^{(m)}_{yx} & \Kb^{(m)}_{yy}
  \end{bmatrix},
  \qquad
  \fb^{(m)} =
  \begin{bmatrix} \fb^{(m)}_{x} & \fb^{(m)}_{y} \end{bmatrix}.
\end{align*}
We introduce the compact notation $\Kb^{(m)}_{\alpha\beta}$ and
$\fb^{(m)}_{\alpha}$ where $\alpha,\beta\in\big\{x,y\big\}$.
Then, we write each component $(\alpha\beta)$ for the global domain
$\Omega$ by considering the contributions from all the $q$ subdomains
in this form
\begin{align}
  \Kb_{\alpha\beta} =
  \begin{bmatrix}
    \Kb_{\alpha\beta,11} & \Kb_{\alpha\beta,12} & \ldots & \Kb_{\alpha\beta,1q} \\
    \Kb_{\alpha\beta,21} & \Kb_{\alpha\beta,22} & \ldots & \Kb_{\alpha\beta,2q} \\
    \vdots & \vdots  & \ddots & \vdots \\
    \Kb_{\alpha\beta,q1} & \Kb_{\alpha\beta,q2} & \ldots & \Kb_{\alpha\beta,qq} 
  \end{bmatrix},
  \label{eq:K}
\end{align}
The block matrices $\Kb_{\alpha\beta,mp}$ for $m$, $p=1,2,\ldots,q$
are originated by the subdomain concatenation through the common nodes
that are shared at the internal interface.
\subsubsection{Subdomain concatenation}
\label{subsubsev:subdomain:concatenation}
Following~\cite{markeeva2021}, a very efficient way to take such a
concatenation into account is to introduce the connectivity matrix
$\Pib^{(mp)}$ for every possible pair of subdomains labeled by $m$ and
$p$, whose $(ij)$-th component is defined as
\begin{align*}
  \big(\Pib^{(mp)}\big)_{ij} =
  \begin{cases}
    1 & \text{if~} i^{m} \sim j^{p}, \\
    0 & \text{otherwise},
  \end{cases}
\end{align*}
where $i^{m}\sim j^{p}$ means that the node $i$ of the domain
$\Omega^{(m)}_h$ coincides with the node $j^{p}$ of the domain
$\Omega^{(p)}_h$.
Let $\gamma$ be a strictly positive, real number.
We compute the diagonal terms of the matrix \eqref{eq:K} as
\begin{align*}
  \Kb_{\alpha\beta,mm} = \Kb^{(m)}_{\alpha\beta} - \gamma \Pib^{(mm)},
\end{align*}
and the extra-diagonal terms ($m\neq p$) as
\begin{subequations}
  \begin{align}            
    \Kb_{\alpha\beta,mp} &= 
    \begin{cases}
      \Pib^{(mp)} \Kb^{(p)}_{\alpha\beta} - \gamma \Pib^{(mp)} & \quad \text{if~} m\text{~and~} p\text{~are~adjacent}, \\
      0 & \quad \text{otherwise}.
    \end{cases}
  \end{align}
\end{subequations}
A convenient value for $\gamma$ is the mean value of the diagonal
element of the stiffness matrix
$\Kb^{\left(m\right)}_{\alpha\beta}$.
Likewise, on employing the connectivity matrices defined above, we
obtain the $\alpha$-component of the global RHS vector $\gb^{(m)}$ for
$\alpha\in\{x,y\}$ starting from the local force vector component
$\fb^{(m)}_\alpha$ associated to the subdomain $\Omega^{(m)}$, which
we write as follows:
\begin{equation}
  \gb^{(m)}_\alpha = \fb^{(m)}_\alpha + \sum_{p\neq m} \Pib^{(mp)}\fb^{(p)}.
\end{equation}
In the above equation, the concatenation term takes into account the
contribution from the other subdomains to the nodes of
$\Omega^{(m)}_h$ that are shared through the common internal
interfaces.

\subsection{Final stiffness matrix and force vector}
According to the previous section, we compute the four global
stiffness matrices $\Kb_{xx}$, $\Kb_{xy}$, $\Kb_{yx}$, $\Kb_{yy}$ and
two global force vectors $\fb_{x}$ and $\fb_{y}$.
Finally, we compute the stiffness matrix $\StfMat$ and the force vector
$\fb$ using the Kronecker product:
\begin{align}
  \StfMat =
  \begin{bmatrix}
    1 & 0\\
    0 & 0
  \end{bmatrix}
  \otimes\Kb_{xx} +
  \begin{bmatrix}
    0 & 1\\
    0 & 0
  \end{bmatrix}
  \otimes\Kb_{xy} +
  \begin{bmatrix}
    0 & 0\\
    1 & 0
  \end{bmatrix}
  \otimes\Kb_{yx} +
  \begin{bmatrix}
    0 & 0\\
    0 & 1
  \end{bmatrix}\otimes\Kb_{yy}
  =
  \begin{bmatrix}
    \Kb_{xx} & \Kb_{xy}\\
    \Kb_{yx} & \Kb_{yy}
  \end{bmatrix},
\end{align}
and
\begin{align}
  \fb =
  \begin{bmatrix}
    1\\
    0
  \end{bmatrix}
  \otimes\fb_{x} +
  \begin{bmatrix}
    0\\
    1
  \end{bmatrix}
  \otimes\fb_{y}
  =
  \begin{bmatrix}
    \fb_{x}\\
    \fb_{y}
  \end{bmatrix}.
\end{align}
Using the stiffness matrix and the RHS load vector defined above, we
solve the resulting linear system $\StfMat\ub=\fb$ for the displacement
field solution $\ub$ on employing the AMEn solver~\cite{dolgov2014}.

\subsection{Dirichlet boundary conditions}
\label{subsec:Dirichlet:boundary:conditions}
We focus on how to impose Dirichlet boundary conditions at the
subdomain level.
Each one of the four sides of every subdomain can be either
constrained via horizontal or vertical rollers or clamped, where both
horizontal and vertical displacements are constrained; see,
cf.~\cite{markeeva2021}.

Overall, the set of all possible combinations of boundary conditions
in the two-dimensional vector case is spanned by the following
one-dimensional vectors:
\begin{align*}
  & \mathbf{X}_{01} = [ 
    \underbrace{0,1,1,\ldots,1,1}_{1,2,3,\ldots,2^d -1, 2^d},\underbrace{1,1,\ldots,1,1}_{1,2,3,\ldots,2^d-1,2^d}
  ],\; \mathbf{X}_{10} = [ 
    \underbrace{1,1,1,\ldots,1,0}_{1,2,3,\ldots,2^d -1, 2^d},\underbrace{1,1,\ldots,1,1}_{1,2,3,\ldots,2^d-1,2^d}
  ],\\
  & \mathbf{Y}_{01} = [ 
    \underbrace{1,1,1,\ldots,1,1}_{1,2,3,\ldots,2^d -1, 2^d},\underbrace{0,1,\ldots,1,1}_{1,2,3,\ldots,2^d-1,2^d}
  ],\; \mathbf{Y}_{10} = [ 
    \underbrace{1,1,1,\ldots,1,1}_{1,2,3,\ldots,2^d -1, 2^d},\underbrace{1,1,\ldots,1,0}_{1,2,3,\ldots,2^d-1,2^d}
  ]\\
  & \mathbf{XY}_{01} = [ 
    \underbrace{0,1,1,\ldots,1,1}_{1,2,3,\ldots,2^d -1, 2^d},\underbrace{0,1,\ldots,1,1}_{1,2,3,\ldots,2^d-1,2^d}
  ],\; \mathbf{XY}_{10} = [ 
    \underbrace{1,1,1,\ldots,1,0}_{1,2,3,\ldots,2^d -1, 2^d},\underbrace{1,1,\ldots,1,0}_{1,2,3,\ldots,2^d-1,2^d}
  ],\\
\end{align*}
where $0$ stands for constrained degrees of freedom (dofs) and $1$ for
free dofs.
When the sides are free, the corresponding boundary vector is
\begin{equation}
  \mathbf{XY}_{11} = [ 
    \underbrace{1,1,1,\ldots,1,1}_{1,2,3,\ldots,2^d -1, 2^d},\underbrace{1,1,\ldots,1,1}_{1,2,3,\ldots,2^d-1,2^d}].
\end{equation}
The boundary masks are then generated through the Kronecker product of
the relevant vectors.
For example, when we eliminate either the nodal $x-$ or
$y-$displacements at the left $(L)$ and right $(R)$ side of the
two-dimensional domain, the masks to be applied are the following:
\begin{align*}
  &\mathbf{M}_{XL} = \mathbf{XY}_{11}^T \otimes \mathbf{X}_{01},\;\mathbf{M}_{XR} = \mathbf{XY}_{11}^T \otimes \mathbf{X}_{10},\\
  &\mathbf{M}_{YL} = \mathbf{XY}_{11}^T \otimes \mathbf{Y}_{01},\;\mathbf{M}_{YR} = \mathbf{XY}_{11}^T \otimes \mathbf{Y}_{10}.
\end{align*}
The z-ordered mask matrices are then obtained by replacing the
Kronecker product $\otimes$ with the z-Kronecker product $\oslash$
\cite{markeeva2021}.
Then, a flatted version of the masks is constructed following the
procedure indicated in \cite{markeeva2021}.

\subsubsection{Algorithm}
We summarize all process algorithmic procedures in
scheme~\ref{algo:1} 

\begin{algorithm}[hbt!]
\caption{Global Finite Element Assembly in \QTT{} format with nodal Z-ordering}
\label{algo:1}
\textbf{Input:}
\begin{itemize}[nosep]
    \item $\Omega$: The computational domain
    \item $q$: Number of subdomains
    \item $d$: Number of \QTT{} dimensions
    \item $\fb(x,y)$, $(x,y)\in\Omega$: forcing term 
    \item $\ub_{0}(x,y)$, $(x,y)\in\partial\Omega$: boundary displacement on $\partial\Omega$ for Dirichlet boundary conditions
\end{itemize}
\textbf{Output:}
\begin{itemize}[nosep]
    \item $\Kb$: Global stiffness matrix
    \item $\fb$: Global force vector
\end{itemize}
\begin{algorithmic}[1]
  \State Partition $\Omega$ into $q$ subdomains $\{\Omega^{(m)}\}_{m=1,2,\ldots,q}$
  \State Generate a local quadrilateral mesh with $2^d\times2^d$ partition elements in each subdomain $\Omega^{(m)}$
  \State Introduce the Z-ordering nodal numbering system in each subdomain
  \For{$m = 1$ \textbf{to} $n$}
  \For{\textbf{each} element $(i,j)$ in $\Omega^{(m)}$} with $i,j=1,2,\ldots,2^d$ 
  \State Compute the Jacobian matrix $\mathbf{J}^{i,j}(\xi,\eta)$ (see the expansion~of Eq.\eqref{eq:Jij:expansion} and subsection~\ref{subsec:subdomain:numerical:integration})
  \State Initialize local stiffness matrix $\Kb^{(m)}$ and force vector $\fb^{(m)}$ to zero.
  \For{ every reference node $(c_1,c_2)\in\{(-1,-1),(1,-1),(1,1),(-1,1)\}$}
  \State Compute matrices $\Kcal_{(c_1,c_2)}$, and $\Gb_{c_1,c_2}$ and shift matrices $\Vb_{c_1}$, $\Vb_{c_2}$
  \State Compute $\StfMat^{(m)}_{(c_1,c_2)}\gets \Vb_{c_1}^T\text{diag}\big(\Kcal_{(c_1,c_2)}\big)\Vb_{c_2},$ (see Eq.~\eqref{eq:Nodal:Stiffness:Matrix})
  \State Compute $\fb^{(m)}_{c_1,c_2}\gets\Vb_{c_1}^T\text{diag}(\Gb_{c_1,c_2})\Vb_{c_2}\overline{\fb}$ (see Eq.~\eqref{eq:load:vector:intensity})
  \State Assemble $\StfMat^{(m)}_{(c_1,c_2)}$ into $\Kb^{(m)}$ (see Eq.~\eqref{eq:Subdomain:Stiffness:Matrix} and subsection~\ref{subsec:assembly:subdomain:level})
  \State Assemble $\fb^{(m)}_{c_1,c_2}$ into $\fb^{(m)}$  (see Eq.~\eqref{eq:subdomain:load:vector:intensity} and subsection~\ref{subsec:assembly:subdomain:level})
  \EndFor
  \EndFor
  \EndFor
  \State Assemble global $\mathbf{K}$ and $\mathbf{f}$ from $\{\mathbf{K}^{(m)}\}$ and $\{\mathbf{f}^{(m)}\}$ (see subsection~\ref{subsec:assembly:global:domain:level})
  \State Apply boundary conditions to $\mathbf{K}$ and $\mathbf{f}$ (see subsection~\ref{subsec:Dirichlet:boundary:conditions})
  \State \Return $\mathbf{K}$, $\mathbf{f}$
\end{algorithmic}
\end{algorithm}

\section{Low-rank TT and QTT formats for solver data representation}
\label{sec4:TT-QTT-formats}

Hereafter, we briefly introduce the concept of tensor
train~\cite{oseledets2011,oseledets2009b} and quantized tensor train
formats~\cite{khoromskij2011,oseledets2010b} for low-rank tensor
representations and discuss how we use these formats in our
algorithm's design.
The key point is that by approximating the PDE solution, which is
represented by a high-dimensional tensor with low-rank TT or QTT
format, we can efficiently capture the essential solution features
without storing the entire solution, hence leading to faster
computations.

\medskip
The \textbf{TT format representation} of a tensor is a compressed
representation that makes it possible to store in memory and
efficiently manipulate high-dimensional low-rank tensors.
This task is achieved by breaking down a $d$-dimensional tensor into a
product of two matrices and $(d-2)$ three-dimensional tensors, called
the \emph{cores}.
Let $\Tcal\in\mathbb{R}^{n_{1}\times\ldots\times n_{d}}$ be a
$d$-dimensional tensor, whose elements are addressed as
$\Tcal(i_1,\ldots,i_d)$, with $i_{\ell}=1,\ldots,n_{\ell}$ for
$\ell=1,\ldots,d$.
We refer to $n_{\ell}$ as the \textit{mode size} (or, simply, the
\textit{size}) of the $\ell$-th dimension, and to the corresponding
index $i_{\ell}$ as the \textit{mode index} (or, simply, the
\textit{index}) associated with such dimension.
Tensor $\Tcal$ is said to be in TT-format if there exists $d$
three-dimensional cores, denoted as
$T_{\ell}\in\mathbb{R}^{r_{\ell-1}\times n_{ell}\times r_{\ell}}$ for
$\ell=1,\ldots,d$, such that
\begin{align*}
  \Tcal(i_{1},i_{2},\ldots,i_{d-1},i_{d}) =
  T_{1}(:,i_{1},:)\,T_{2}(:,i_{2},:)\ldots T_{d-1}(:,i_{d-1},:)\,T_{d}(:,i_{d},:)
\end{align*}
for all possible combinations of the indices $i_1,i_2,\ldots,i_d$, and
using a Matlab-like notation to express the matrix-matrix
multiplication.
The numbers $r_{\ell}$ for $\ell=0,1,\ldots,d$ are called the
\textit{representation ranks} (or, simply, the \textit{ranks}) of the
TT-format representation.
We assume that $r_0=r_d=1$, so that the first and last core, i.e.,
$T_1$ and $T_d$ are, indeed, matrices.
By expanding the compact core-matrix representation above (and
omitting the trivial summation over the indices $\alpha_0$ and
$\alpha_d$ ranging from $1$ to $1$), we find the equivalent entry-wise
formulation:
\begin{multline*}
  \qquad
  \Tcal(i_{1},i_{2}\ldots i_{d-1},i_{d}) =
  \sum_{\alpha_1=1}^{r_1}\sum_{\alpha_2=1}^{r_2}\ldots\sum_{\alpha_{d-1}=1}^{r_{d-1}}
  T_{1}(i_{1},\alpha_1)\,T_{2}(\alpha_1,i_{2},\alpha_2)\ldots\\
  \ldots T_{d-1}(\alpha_{d-2},i_{d-1},\alpha_{d-1})\,T_{d}(\alpha_{d-1},i_{d}).
  \qquad
\end{multline*}
The tensor train decomposition is particularly useful in the case
where the ranks are much smaller than the mode sizes, i.e.,
$r_{\ell}\ll n_{\ell}$ for all tensor dimensions $\ell=1,\ldots,d$.
The amount of memory required for storing $T$ is
$\textsf{storage}=r_{1}n_{1}+\sum_{i=2}^{d-1}r_{i-1}n_{i}r_{i}+r_{d-1}
n_{d}$, which grows proportionally to $\mathcal{O}(dr^2n)$, where
$r=\max(r_1,\ldots,r_{d-1})$ and $n=\max(n_1,\ldots,n_{d})$ are
convenient upper bounds on the ranks and mode sizes.
Following~\cite{markeeva2021}, we introduce the
\textit{effective~rank}, which we denote by $r_e$ and occasionally
abbreviate as \textit{erank}.
This quantity is a sort of average rank for the cores of the tensor
train representation of $\Tcal$ and satisfies, by definition, the
condition
\begin{align*}
  r_{e} n_{1}+\sum_{i=2}^{d-1} r_{e}^{2} n_{i}+r_{e} n_{d} = \textsf{storage of $T$}.
\end{align*}
On solving the equation above, it follows immediately that
$r_e=\mathcal{O}(\sqrt{\textsf{storage of $T$}})$

In the case of second-order elliptic partial differential equations in
two dimensions \cite{kazeev2015}, the quantized tensor-structured
solver applied to $q$ subdomains has been characterized in terms of
the maximum rank $R_{d}$
\begin{equation}\label{eq:Rd}
    R_{d}=\max_{0\leq i < d} r_{d},
\end{equation}
and the number of parameters $N_{d}$ involved in the representation defined as follows:
\begin{equation}\label{eq:Nd}
    N_{d}=q r_{0}+\sum^{d-1}_{i=1} 2^{2}r_{i-1}r_{i}+2^{2}r_{d-1}=\mathcal{O}\left(d R^{2}_{d}\right).
\end{equation}

Now, consider the \textit{multi-index representation} of the tensor
elements.
Let $\iu=(i_1,\ldots,i_d)$ a $d$-dimensional multi-index, where every
component $i_{\ell}$ is a non-negative integer number, and
$\abs{\iu}=i_1+\ldots+i_d$ is the \textit{order} of $\ib$.
Then, we address the elements of tensor $T$ as
$T(\iu):=T(i_1,\ldots,T_d)$.
Since $T$ has a single multi-index, we can interpret it as a
``\textit{vector}'' tensor (a vector is an object with a single
index), and it is straightforward to introduce the concept of a
``matrix'' tensor, i.e., a \textit{multi-dimensional matrix} whose
elements are indexed by two multi-indices $\iu$ and $\ju$, e.g.,
$\Acal(\iu,\ju)$.
For consistency, these multi-indices must have the same number of
dimensions $d$, so that we can pair them as $(i_{\ell},j_{\ell})$ for
$\ell=1,\ldots,d$, although the mode sizes can be different.
Let $i_{\ell}=1,\ldots,n_{\ell}$ (as above), and
$j_{\ell}=1,\ldots,m_{\ell}$ for $\ell=1,\ldots,m_{\ell}$.
A $d$-dimensional matrix tensor $\Acal$ is said to be in tensor-train
format if there exist $d$ four-dimensional tensor said \textit{cores}
$A_{\ell}\in\mathbb{R}^{r_{\ell-1},n_{\ell},m_{\ell},r_{\ell}}$ with
$r_1=r_d=1$ such that
\begin{align*}
  \Acal(\iu,\ju)
  & := \Acal\big((i_1,\ldots,i_d),(j_1,\ldots,j_d)\big) \\[0.5em]
  &=
  A_{1}(:,i_1,j_1,:)A_{2}(:,i_2,j_2,:) \cdots 
  A_{d-1}(:,i_{d-1},j_{d-1},:)A_{d}(:,i_{d},j_{d},:)
\end{align*}
for all possible combinations of the index pairs
$(i_1,j_1),\ldots,(i_d,j_d)$, using again a Matlab-like notation to
express the matrix-matrix multiplication.
The equivalent entry-wise expression is
\begin{multline*}
  \qquad
  \Acal\big((i_1,\ldots,i_d),(j_1,\ldots,j_d)\big) = \\ =
  \sum_{\alpha=1}^{r_1}\sum_{\alpha_2=1}^{r_2}\ldots\sum_{\alpha_{d-2}=1}^{r_{d-2}}\sum_{\alpha_{d-1}=1}^{r_{d-1}}
  A_{1}(i_1,j_1,\alpha_1)A_{2}(\alpha_1,i_2,j_2,\alpha_2)
  \cdots \\ \cdots 
  A_{d-1}(\alpha_{d-2},i_{d-1},j_{d-1},\alpha_{d-1})A_{d}(\alpha_{d-1},i_{d},j_{d}),
  \qquad
\end{multline*}
and again omitting the trivial summation over $\alpha_0$ and
$\alpha_d$.

\begin{remark}
  We can treat each pair $i_{\ell},j_{\ell}$ as one ``long index'',
  e.g., $k_{\ell}=i_{\ell} + n_{\ell}(j_\ell-1)$ ranging from $1$ to
  $n_{\ell}m_{\ell}$.
  Moreover, we can consider each core $A_{\ell}$ as a block matrix of
  dimensions $(r_{\ell-1}n_{\ell})\times(r_{\ell}m_{\ell})$.
  This block matrix is split into an $r_{\ell-1}\times r_{\ell}$ block
  structure.
  The block $A_{\ell}(\alpha_{\ell-1},:,:,\alpha_{\ell})$, which is
  addressed by the index pair $(\alpha_{\ell-1},\alpha_{\ell})$, for
  $\alpha_{\ell-1}=1,2,\ldots,r_{\ell-1}$,
  $\alpha_{\ell}=1,2,\ldots,r_{\ell}$, is the $n_{\ell}\times
  m_{\ell}$-sized matrix block, whose $(i_{\ell},j_{\ell})$-th element
  is $A_{\ell}(\alpha_{\ell-1},i_{\ell},j_{\ell},\alpha_{\ell})$.
\end{remark}

An alternative formulation makes use of the \emph{strong Kronecker
product} denoted as $\bowtie$.
The strong Kronecker product is a block matrix operation applied to
block sub-matrices, similar to the regular matrix product, but with
the key difference that it multiplies entire blocks using the
(regular) Kronecker product $\otimes$ instead of multiplying
individual elements.
For example, consider matrices $K=\big(K_{IJ}\big)$ and
$L=\big(L_{I'J'}\big)$ that are partitioned into smaller sub-matrices
(or blocks) $K_{Ij}$ and $L_{I'J'}$ whose position inside,
respectively, $K$ and $L$ is determined by the block index pair ``$IJ$''
and ``$I'J'$''.
In the strong Kronecker product between matrices $K$ and $L$ the
element at position $\big( (IJ), (I'J') \big)$ is the Kronecker
product of the corresponding subblocks, and is formally expressed as:
\begin{align*}
  (K \bowtie L)_{ \big( (IJ), (I'J') \big)} = K_{IJ} \otimes L_{I'J'}.
\end{align*}
Using the strong Kronecker product, we can write $\Acal$ in TT format
in terms of its cores $\{A_{\ell}\}_{1\leq\ell\leq d}$ as
\begin{align*}
  \Acal = A_{1}\bowtie A_{2}\bowtie\ldots\bowtie A_{d-1}\bowtie A_{d}.
\end{align*}

\begin{remark}
  If all the internal ranks of tensor $\Acal$ are equal to $1$, we can
  write it as the regular Kronecker product of the $d$ matrices
  $A_{\ell}(1,:,:,1)\in\mathbb{R}^{n_{\ell}\times m_{\ell}}$,
  $\ell=1,2,\ldots,d$:
  \begin{align*}
    \Acal = A_{1} \otimes A_{2} \otimes \ldots \otimes A_{d-1} \otimes A_{d},
  \end{align*}
  with the obvious extension of the Kronecker product from two to $d$
  arguments.
\end{remark}

\medskip
Following~\cite{markeeva2021}, we will use the \textbf{QTT-format
  representation} and adopt the ``\textit{Z-order}'' that we will
extend to the vector case in the next section.
QTT is a special tensor train format that can be applied to represent
both vector and matrices.
To obtain a QTT representation of such objects, we first reshape them
into a multidimensional \textit{binary} tensor representation, where
``binary'' that the size of every dimension is $2$, and then we
perform a tensor train decomposition.
The motivation for using QTT is its low memory consumption, with
low-rank approximations and speed in solving linear equation
systems.
We will employ the QTT format to compute and store all the
components of the final linear system.
Then, we will solve this system by using the AMEn solver, which is
designed explicitly for TT representations to take
advantage of the low-rank structure when present, hence leading to
faster convergence and reduced computation time.

Using QTT improves efficiency in solving differential equations
through several key mechanisms.
QTT reduces the complexity of representing matrix operators and
solution vectors from polynomial, e.g., $\mathcal{O}(N^2)$, to
logarithmic, e.g., $\mathcal{O}(\log(N))$, where $N$ is the total
number of degrees of freedom.
Such reduction to a logarithmic complexity drastically decreases both
memory usage and computational cost.
For example, we can reshape a vector $\vb\in\mathbb{R}^{2^{d}}$ as a
$d$-dimensional tensor $\Vcal$ of dimension
$2\times2\times\ldots\times2$ ($d$ times); we consider its TT format
representation with ranks $r_1,\ldots,r_{d-1}$; we reshape vector
$\vb$ into the $d$-dimensional tensor $\Vcal$ encoding the index value
of $\vb$, e.g., $1\leq i\leq 2^{d}$, into the binary format:
\begin{align*}
  i
  = \overline{i_{1}, i_{2}, \ldots i_{d}}
  = \sum_{k=1}^{d} 2^{k-1} i_{k}
  \leftrightarrow \left(i_{1}, i_{2}, \ldots i_{d}\right)
\end{align*}
We can use the same idea to represent a matrix in the QTT format.
In such a case, we formally find that
$\Acal\in\mathbb{R}^{2^{d}\times2^{d}}$.
According to \cite{markeeva2021}, we will reorder the matrices
elements in the \textit{Z-order} to prevent unnecessary rank growth
during computations and make use of operations like \textit{z-kron} to
manipulate these mathematical objects efficiently.

\section{Numerical examples}
\label{sec5:numerical}

\newcommand{\FIGREF}[2]{\ref{#1}$(#2)$}
\newcommand{\FEniCS}{\textrm{FEniCS}}

\begin{figure}[htb!]
  \centering
  \begin{tabular}{c}
    \includegraphics[width=1\textwidth]{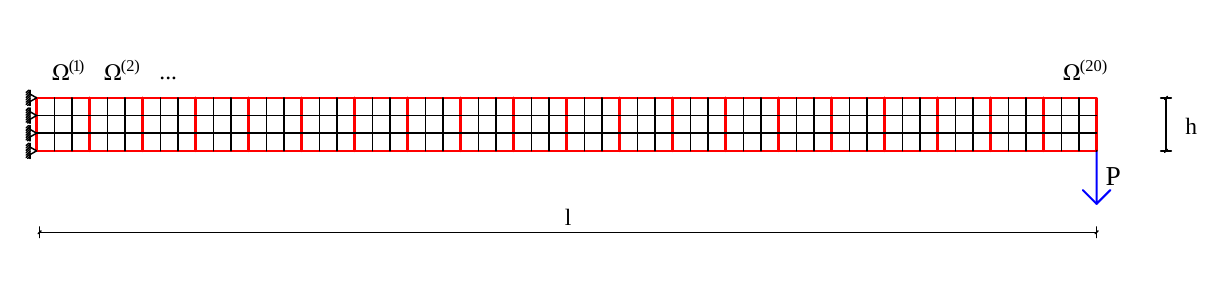}\\[-1.em]
    \hspace{-1cm}$\mathbf{(a)}$
  \end{tabular}
  \begin{tabular}{cc}
    \includegraphics[width=0.5\textwidth,trim={0cm 0cm 0 0cm}]{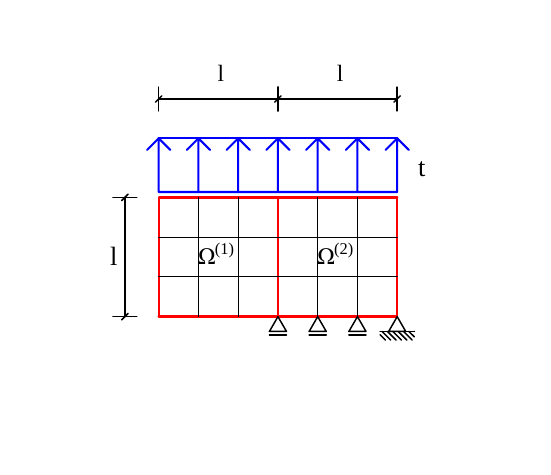} & \hspace{-2cm}
    \includegraphics[width=0.5\textwidth, trim={0 0cm 0 0}]{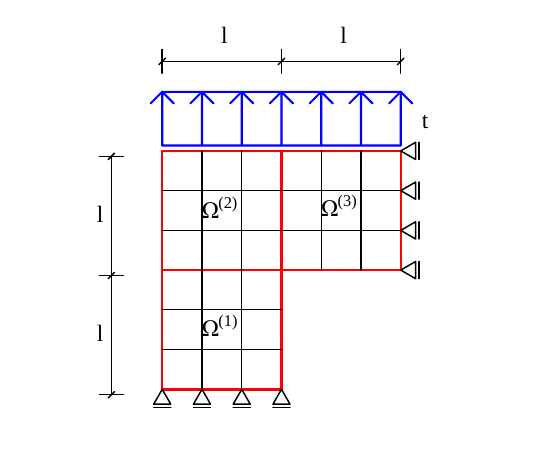} \\[-1.5em]
    $\mathbf{(b)}$ & \hspace{-2cm} $\mathbf{(c)}$
  \end{tabular}
  \caption{Geometry, boundary conditions, and representative meshes
    for the test cases of Section~\ref{sec5:numerical}.
    $(a)$ Cantilever beam test case assuming $d=2$ levels and $q=20$
    subdomains.
    $(b)$ Single edge notch tensile test case with $d=2$ levels and
    $q=2$ subdomains.
    $(c)$ L-shaped domain test case with $d=3$ levels and $q=3$
    subdomains.
    Red vertical and horizontal lines separate adjacent subdomains;
    each subdomain is partitioned by a $3\times3$-square grid, and
    denoted by $\Omega^{(m)}$.
    }
  \label{fig:geometries}
\end{figure}

This section is devoted to the assessment of the convergence
properties of QTT-FEM compared to classic sparse FEM solvers in a
series of structural plane stress examples that require accurate
discretization to overcome the poor performance of low-order elements
such as in beam bending problems or because of geometric
singularities.
In the former case, bending in the cantilever beam as in
Figure~\FIGREF{fig:geometries}{a} is studied having a height-to-length
ratio equal to $20$.
For this geometry, 20 subdomains were used, each with a grid of
$2^d\times 2^d$ nodes.
A red vertical line separates adjacent subdomains, and each subdomain
is partitioned by a $3\times3$-square grid and denoted by
$\Omega^{(m)}$, $m = 1,2,\ldots,20$.
Then, a standard tensile single-edge notched (SEN) specimen and an
L-shaped plate are modeled as examples of possible applications
relevant to solid mechanics engineering, see,
e.g.,~\cite{zienkiewicz2005}.
The corresponding geometries are displayed in Figures
\FIGREF{fig:geometries}{b}-\FIGREF{fig:geometries}{c}, assuming $d=2$
levels and $q=2$ subdomains and $d=3$ levels and $3$ subdomains,
respectively.
Both the SEN and the L-shaped plate have a length of $\ell=1$ mm and
are subjected to a traction $t=3$ MPa.
In all these examples, we assume the material to be isotropic with
Young modulus E = 64 MPa and a vanishing Poisson ratio.
  
\subsection{Convergence properties}
 
The present section is dedicated to the study of the convergence
properties of the QTT finite element solver.
In particular, we compare the performance of our QTT approach with
that of \FEniCS{}~\cite{baratta2023}, a classical sparse matrix finite
element software.

According to~\cite{zienkiewicz2005,bathe2006,szabo2021}, we evaluate
$E=\|\ub_{QTT}-\ub_{\refsol}\|$, the error in the energy seminorm, and
$E_{\mathbb{L}^2}=\|\ub_{QTT}-\ub_{\refsol{}}\|_{\mathbb{L}^2},$, the
error in the $\mathbb{L}^2$-norm.
The notation $\ub_{QTT}$ identifies the displacement solution obtained
with the QTT solver, while the reference solution $\ub_{\refsol}$ is
the displacement evaluated with \FEniCS{} with an overrefined mesh
consisting of $7242852$ degrees of freedom (dofs).
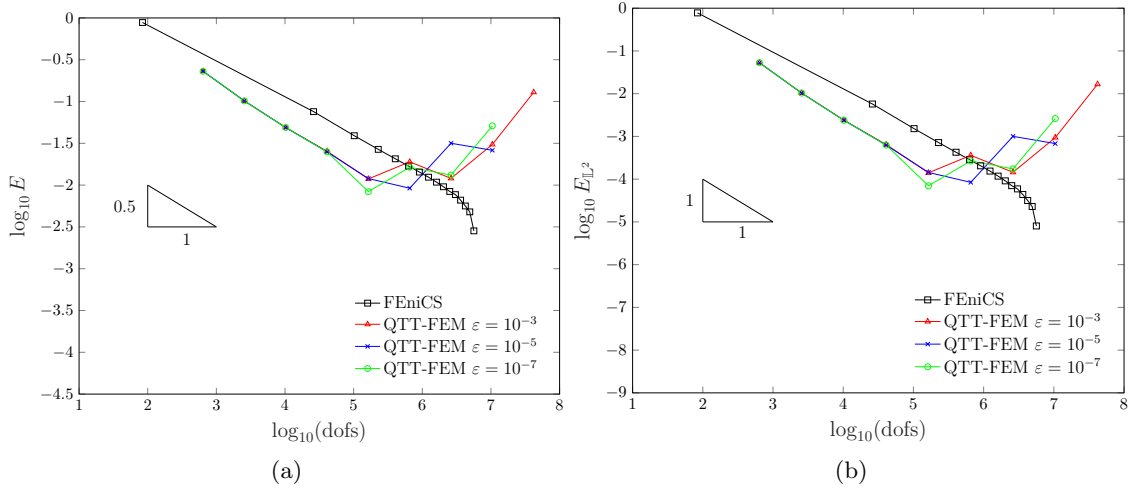
\begin{figure}[htb!]
  \centering
  \subfloat[]{\resizebox{.45\linewidth}{!}{
%
%
\begin{tikzpicture}

\begin{axis}[%
width=4.521in,
height=3.566in,
at={(0.758in,0.481in)},
scale only axis,
unbounded coords=jump,
xmin=1,
xmax=8,
xlabel style={font=\color{white!15!black}},
xlabel={\Large $\log_{10}$(dofs)},
ymin=-4.5,
ymax=0,
ylabel style={font=\color{white!15!black}},
ylabel={\Large $\log_{10} E$},
axis background/.style={fill=white},
legend style={at={(0.97,0.03)}, anchor=south east, legend cell align=left, align=left, fill=none, draw=none}
]
\addplot [color=black, mark=square, mark options={solid, black}]
  table[row sep=crcr]{%
1.92427928606188	-0.0547875356247404\\
4.41584106950205	-1.1207909577374\\
5.00903424924124	-1.40909065897324\\
5.35822428018804	-1.57214793916155\\
5.60659845808283	-1.68513680464531\\
5.79951423251205	-1.77200908796785\\
5.957272938382	-1.84379568461069\\
6.09073477424316	-1.90531439177495\\
6.20639460446347	-1.96460214937915\\
6.30844744774986	-2.01982780343447\\
6.39976057169269	-2.07590561432008\\
6.48238072267896	-2.11384107463041\\
6.55782011802392	-2.1797395122972\\
6.6272277617714	-2.24958205752615\\
6.69149718978175	-2.3201458925503\\
6.75133700605421	-2.54731800743226\\
};
\addlegendentry{\large FEniCS}

\addplot [color=red, mark=triangle, mark options={solid, red}]
  table[row sep=crcr]{%
2.80617997398389	-0.638137314860492\\
3.40823996531185	-0.992580841394672\\
4.01029995663981	-1.31044953678371\\
4.61235994796777	-1.59538277594145\\
5.21441993929574	-1.92576429383878\\
5.8164799306237	-1.72160468873907\\
6.41853992195166	-1.91931390449031\\
7.02059991327962	-1.51359395946414\\
7.62265990460759	-0.889438753292321\\
};
\addlegendentry{\large QTT-FEM $\varepsilon=10^{-3}$}

\addplot [color=blue, mark=x, mark options={solid, blue}]
  table[row sep=crcr]{%
2.80617997398389	-0.637827571540769\\
3.40823996531185	-0.992756798746244\\
4.01029995663981	-1.31129149608763\\
4.61235994796777	-1.60174218650446\\
5.21441993929574	-1.92134785265516\\
5.8164799306237	-2.0372162000892\\
6.41853992195166	-1.49917828669313\\
7.02059991327962	-1.5836102510518\\
};
\addlegendentry{\large QTT-FEM $\varepsilon=10^{-5}$}

\addplot [color=green, mark=o, mark options={solid, green}]
  table[row sep=crcr]{%
2.80617997398389	-0.637824719324442\\
3.40823996531185	-0.992749089072084\\
4.01029995663981	-1.31115704487425\\
4.61235994796777	-1.60147639868783\\
5.21441993929574	-2.07783829450049\\
5.8164799306237	-1.78622131049067\\
6.41853992195166	-1.87970318523689\\
7.02059991327962	-1.29053387758547\\
};
\addlegendentry{\large QTT-FEM $\varepsilon=10^{-7}$}

\addplot [color=black, forget plot]
  table[row sep=crcr]{%
2	-2\\
2.30102999566398	-2.15051499783199\\
2.47712125471966	-2.23856062735983\\
2.60205999132796	-2.30102999566398\\
2.69897000433602	-2.34948500216801\\
2.77815125038364	-2.38907562519182\\
2.84509804001426	-2.42254902000713\\
2.90308998699194	-2.45154499349597\\
2.95424250943932	-2.47712125471966\\
3	-2.5\\
};
\addplot [color=black, forget plot]
  table[row sep=crcr]{%
2	-2\\
2	-2.15051499783199\\
2	-2.23856062735983\\
2	-2.30102999566398\\
2	-2.34948500216801\\
2	-2.38907562519182\\
2	-2.42254902000713\\
2	-2.45154499349597\\
2	-2.47712125471966\\
2	-2.5\\
};
\addplot [color=black, forget plot]
  table[row sep=crcr]{%
2	-2.5\\
2.30102999566398	-2.5\\
2.47712125471966	-2.5\\
2.60205999132796	-2.5\\
2.69897000433602	-2.5\\
2.77815125038364	-2.5\\
2.84509804001426	-2.5\\
2.90308998699194	-2.5\\
2.95424250943932	-2.5\\
3	-2.5\\
};
\node[right, align=left, inner sep=0]
at (axis cs:1.5,-2.25) {\large 0.5};
\node[right, align=left, inner sep=0]
at (axis cs:2.5,-2.65) {\large 1};

\end{axis}
\end{tikzpicture}%
}}
  \subfloat[]{\resizebox{.45\linewidth}{!}{
%
%
\begin{tikzpicture}

\begin{axis}[%
width=4.521in,
height=3.566in,
at={(0.758in,0.481in)},
scale only axis,
unbounded coords=jump,
xmin=1,
xmax=8,
xlabel style={font=\color{white!15!black}},
xlabel={\Large $\log_{10}$(dofs)},
ymin=-9,
ymax=0,
ylabel style={font=\color{white!15!black}},
ylabel={\Large $\log_{10}E_{\mathbb{L}^2}$},
axis background/.style={fill=white},
legend style={at={(0.97,0.03)}, anchor=south east, legend cell align=left, align=left, fill=none, draw=none}
]
\addplot [color=black, mark=square, mark options={solid, black}]
  table[row sep=crcr]{%
1.92427928606188	-0.109575071249481\\
4.41584106950205	-2.2415819154748\\
5.00903424924124	-2.81818131794648\\
5.35822428018804	-3.1442958783231\\
5.60659845808283	-3.37027360929062\\
5.79951423251205	-3.5440181759357\\
5.957272938382	-3.68759136922139\\
6.09073477424316	-3.8106287835499\\
6.20639460446347	-3.92920429875831\\
6.30844744774986	-4.03965560686893\\
6.39976057169269	-4.15181122864016\\
6.48238072267896	-4.22768214926082\\
6.55782011802392	-4.35947902459441\\
6.6272277617714	-4.49916411505229\\
6.69149718978175	-4.64029178510061\\
6.75133700605421	-5.09463601486452\\
};
\addlegendentry{\large FEniCS}

\addplot [color=red, mark=triangle, mark options={solid, red}]
  table[row sep=crcr]{%
2.80617997398389	-1.27627462972098\\
3.40823996531185	-1.98516168278934\\
4.01029995663981	-2.62089907356742\\
4.61235994796777	-3.19076555188289\\
5.21441993929574	-3.85152858767756\\
5.8164799306237	-3.44320937747815\\
6.41853992195166	-3.83862780898061\\
7.02059991327962	-3.02718791892828\\
7.62265990460759	-1.77887750658464\\
};
\addlegendentry{\large QTT-FEM $\varepsilon=10^{-3}$}

\addplot [color=blue, mark=x, mark options={solid, blue}]
  table[row sep=crcr]{%
2.80617997398389	-1.27565514308154\\
3.40823996531185	-1.98551359749249\\
4.01029995663981	-2.62258299217525\\
4.61235994796777	-3.20348437300891\\
5.21441993929574	-3.84269570531032\\
5.8164799306237	-4.07443240017839\\
6.41853992195166	-2.99835657338626\\
7.02059991327962	-3.16722050210361\\
};
\addlegendentry{\large QTT-FEM $\varepsilon=10^{-5}$}

\addplot [color=green, mark=o, mark options={solid, green}]
  table[row sep=crcr]{%
2.80617997398389	-1.27564943864888\\
3.40823996531185	-1.98549817814417\\
4.01029995663981	-2.62231408974849\\
4.61235994796777	-3.20295279737566\\
5.21441993929574	-4.15567658900097\\
5.8164799306237	-3.57244262098134\\
6.41853992195166	-3.75940637047378\\
7.02059991327962	-2.58106775517093\\
};
\addlegendentry{\large QTT-FEM $\varepsilon=10^{-7}$}

\addplot [color=black, forget plot]
  table[row sep=crcr]{%
2	-4\\
2.30102999566398	-4.30102999566398\\
2.47712125471966	-4.47712125471966\\
2.60205999132796	-4.60205999132796\\
2.69897000433602	-4.69897000433602\\
2.77815125038364	-4.77815125038364\\
2.84509804001426	-4.84509804001426\\
2.90308998699194	-4.90308998699194\\
2.95424250943932	-4.95424250943933\\
3	-5\\
};
\addplot [color=black, forget plot]
  table[row sep=crcr]{%
2	-4\\
2	-4.30102999566398\\
2	-4.47712125471966\\
2	-4.60205999132796\\
2	-4.69897000433602\\
2	-4.77815125038364\\
2	-4.84509804001426\\
2	-4.90308998699194\\
2	-4.95424250943933\\
2	-5\\
};
\addplot [color=black, forget plot]
  table[row sep=crcr]{%
2	-5\\
2.30102999566398	-5\\
2.47712125471966	-5\\
2.60205999132796	-5\\
2.69897000433602	-5\\
2.77815125038364	-5\\
2.84509804001426	-5\\
2.90308998699194	-5\\
2.95424250943932	-5\\
3	-5\\
};
\node[right, align=left, inner sep=0]
at (axis cs:1.75,-4.5) {\large 1};
\node[right, align=left, inner sep=0]
at (axis cs:2.5,-5.25) {\large 1};

\end{axis}
\end{tikzpicture}%
}}
  \caption{Cantilever beam test case: bilogarithmic plot of QTT and
    \FEniCS{} energy seminorm $(a)$ and $\mathbb{L}^{2}$-norm $(b)$
    errors for different AMEn approximation accuracy.}
  \label{fig:energy_error_cantilever}
\end{figure}
Figures
\FIGREF{fig:energy_error_cantilever}{a}-\FIGREF{fig:energy_error_cantilever}{b}
illustrate the convergence rates in the cantilever test for different
levels of approximation tolerance $\varepsilon$ of the AMEn solver.
In these figures, we plot the errors in both the energy seminorm and
the $\mathbb{L}^{2}$-norm.
For comparison, we show the results obtained by using \FEniCS and the
expected slopes of convergence for the P1-approximation.
The convergence slope for the L-shaped panel concurs with the
literature value \cite{zienkiewicz2005}.
The trend of the QTT solver results is consistent with the expected
convergence rates up to a certain value of the number of degrees of freedom beyond which it drifts
apart because of the rounding errors stemming from both limited machine
precision and AMEn approximation
prevail~\cite{markeeva2021,bachmayr2020}.
  
\begin{figure}[htb!]
  \centering
  \subfloat[]{\resizebox{.45\linewidth}{!}{
%
%
\begin{tikzpicture}

\begin{axis}[%
width=4.521in,
height=3.566in,
at={(0.758in,0.481in)},
scale only axis,
unbounded coords=jump,
xmin=1,
xmax=8,
xlabel style={font=\color{white!15!black}},
xlabel={\Large $\log_{10}$(dofs)},
ymin=-3,
ymax=0,
ylabel style={font=\color{white!15!black}},
ylabel={\Large $\log_{10} E$},
axis background/.style={fill=white},
legend style={legend cell align=left, align=left, fill=none, draw=none}
]
\addplot [color=black, mark=square, mark options={solid, black}]
  table[row sep=crcr]{%
1.7481880270062	-0.224195456750984\\
3.82229887126237	-0.737191996098271\\
4.1691452009912	-0.827212271498564\\
4.41634089057448	-0.892219208180146\\
4.60854742686711	-0.943458074790881\\
5.20737042570465	-1.1094803963619\\
5.8078068843684	-1.29598185736126\\
6.15944747506963	-1.42750666318197\\
6.40905384366752	-1.54450006649999\\
6.60271116177951	-1.66598877436097\\
6.76096516312391	-1.81648436321794\\
6.89478123977332	-2.08848815411733\\
6.9546766833321	-inf\\
};
\addlegendentry{\large FEniCS}

\addplot [color=red, mark=triangle, mark options={solid, red}]
  table[row sep=crcr]{%
1.80617997398389	-0.316407123259727\\
2.40823996531185	-0.480643372137515\\
3.01029995663981	-0.640671288954301\\
3.61235994796777	-0.799284725778552\\
4.21441993929574	-0.960060218371031\\
4.8164799306237	-1.12906638358868\\
5.41853992195166	-1.32030705488656\\
6.02059991327962	-1.5853089118707\\
6.62265990460759	-1.93912353085003\\
7.22471989593555	-1.63863156746099\\
};
\addlegendentry{\large QTT-FEM $\varepsilon=10^{-3}$}

\addplot [color=blue, mark=x, mark options={solid, blue}]
  table[row sep=crcr]{%
1.80617997398389	-0.316645655968644\\
2.40823996531185	-0.480613167941138\\
3.01029995663981	-0.640763119293502\\
3.61235994796777	-0.799293972699999\\
4.21441993929574	-0.960141514081948\\
4.8164799306237	-1.12935438290501\\
5.41853992195166	-1.32077567970019\\
6.02059991327962	-1.58579229492837\\
6.62265990460759	-1.94369280319639\\
7.22471989593555	-1.63804171894876\\
};
\addlegendentry{\large QTT-FEM $\varepsilon=10^{-5}$}

\addplot [color=green, mark=o, mark options={solid, green}]
  table[row sep=crcr]{%
1.80617997398389	-0.316644306566964\\
2.40823996531185	-0.480614495039412\\
3.01029995663981	-0.640764378739109\\
3.61235994796777	-0.79929346855353\\
4.21441993929574	-0.960140707096115\\
4.8164799306237	-1.12935225985424\\
5.41853992195166	-1.32078271018572\\
6.02059991327962	-1.58579924365566\\
6.62265990460759	-1.94390702496207\\
7.22471989593555	-1.63798995454493\\
};
\addlegendentry{\large QTT-FEM $\varepsilon=10^{-7}$}

\addplot [color=black, forget plot]
  table[row sep=crcr]{%
2	-1.5\\
2.30102999566398	-1.575257498916\\
2.47712125471966	-1.61928031367992\\
2.60205999132796	-1.65051499783199\\
2.69897000433602	-1.674742501084\\
2.77815125038364	-1.69453781259591\\
2.84509804001426	-1.71127451000356\\
2.90308998699194	-1.72577249674799\\
2.95424250943932	-1.73856062735983\\
3	-1.75\\
};
\addplot [color=black, forget plot]
  table[row sep=crcr]{%
2	-1.5\\
2	-1.575257498916\\
2	-1.61928031367992\\
2	-1.65051499783199\\
2	-1.674742501084\\
2	-1.69453781259591\\
2	-1.71127451000356\\
2	-1.72577249674799\\
2	-1.73856062735983\\
2	-1.75\\
};
\addplot [color=black, forget plot]
  table[row sep=crcr]{%
2	-1.75\\
2.30102999566398	-1.75\\
2.47712125471966	-1.75\\
2.60205999132796	-1.75\\
2.69897000433602	-1.75\\
2.77815125038364	-1.75\\
2.84509804001426	-1.75\\
2.90308998699194	-1.75\\
2.95424250943932	-1.75\\
3	-1.75\\
};
\node[right, align=left, inner sep=0]
at (axis cs:1.45,-1.6) {\large 0.25};
\node[right, align=left, inner sep=0]
at (axis cs:2.5,-1.9) {\large 1};

\end{axis}
\end{tikzpicture}%
}}
  \subfloat[]{\resizebox{.45\linewidth}{!}{
%
%
\begin{tikzpicture}

\begin{axis}[%
width=4.521in,
height=3.566in,
at={(0.758in,0.481in)},
scale only axis,
unbounded coords=jump,
xmin=1,
xmax=8,
xlabel style={font=\color{white!15!black}},
xlabel={\Large $\log_{10}$(dofs)},
ymin=-3,
ymax=0,
ylabel style={font=\color{white!15!black}},
ylabel={\Large $\log_{10} E$},
axis background/.style={fill=white},
legend style={legend cell align=left, align=left, fill=none, draw=none}
]
\addplot [color=black, mark=square, mark options={solid, black}]
  table[row sep=crcr]{%
1.90308998699194	-0.391196209924492\\
3.40857912540867	-1.02143863646458\\
3.75143308181935	-1.15808101521747\\
3.99659922270017	-1.25264557943798\\
4.18757711905509	-1.32464844689989\\
4.78391786504047	-1.54246917616074\\
5.38310052508721	-1.75941765011451\\
5.73432128223689	-1.89458063665146\\
5.98371737568707	-2.00120906571443\\
6.1772484137745	-2.0964582714189\\
6.33541818048625	-2.18974645825851\\
6.46917406603659	-2.28973040497331\\
6.58505467180525	-2.40996662925559\\
6.6872793708305	-2.5887133873178\\
6.77873006856815	-inf\\
};
\addlegendentry{\large FEniCS}

\addplot [color=red, mark=triangle, mark options={solid, red}]
  table[row sep=crcr]{%
1.98227123303957	-0.594012652224472\\
2.58433122436753	-0.874350241572389\\
3.18639121569549	-1.12126778002076\\
3.78845120702346	-1.34189518776249\\
4.39051119835142	-1.54954685018664\\
4.99257118967938	-1.75633205484794\\
5.59463118100734	-2.06576094381956\\
6.19669117233531	-2.71509532210441\\
6.79875116366327	-2.3096424906186\\
7.40081115499123	-2.16460623031882\\
};
\addlegendentry{\large QTT-FEM $\varepsilon=10^{-3}$}

\addplot [color=blue, mark=x, mark options={solid, blue}]
  table[row sep=crcr]{%
1.98227123303957	-0.592496064934009\\
2.58433122436753	-0.874850200321086\\
3.18639121569549	-1.12209052408798\\
3.78845120702346	-1.34535733547249\\
4.39051119835142	-1.55511268505851\\
4.99257118967938	-1.7641634750365\\
5.59463118100734	-1.99799411989371\\
6.19669117233531	-2.38743538917541\\
6.79875116366327	-2.34457207119579\\
7.40081115499123	-2.21156024306446\\
};
\addlegendentry{\large QTT-FEM $\varepsilon=10^{-5}$}

\addplot [color=green, mark=o, mark options={solid, green}]
  table[row sep=crcr]{%
1.98227123303957	-0.592489780852644\\
2.58433122436753	-0.874853001605566\\
3.18639121569549	-1.12209091572335\\
3.78845120702346	-1.3453819135876\\
4.39051119835142	-1.55516031907771\\
4.99257118967938	-1.76417686928286\\
5.59463118100734	-1.99827934606971\\
6.19669117233531	-2.38531910869114\\
6.79875116366327	-2.34619041607711\\
7.40081115499123	-2.21219848930615\\
};
\addlegendentry{\large QTT-FEM $\varepsilon=10^{-7}$}

\addplot [color=black, forget plot]
  table[row sep=crcr]{%
2	-1.9\\
2.30102999566398	-2.03546349804879\\
2.47712125471966	-2.11470456462385\\
2.60205999132796	-2.17092699609758\\
2.69897000433602	-2.21453650195121\\
2.77815125038364	-2.25016806267264\\
2.84509804001426	-2.28029411800642\\
2.90308998699194	-2.30639049414637\\
2.95424250943932	-2.3294091292477\\
3	-2.35\\
};
\addplot [color=black, forget plot]
  table[row sep=crcr]{%
2	-1.9\\
2	-2.03546349804879\\
2	-2.11470456462385\\
2	-2.17092699609758\\
2	-2.21453650195121\\
2	-2.25016806267264\\
2	-2.28029411800642\\
2	-2.30639049414637\\
2	-2.3294091292477\\
2	-2.35\\
};
\addplot [color=black, forget plot]
  table[row sep=crcr]{%
2	-2.35\\
2.30102999566398	-2.35\\
2.47712125471966	-2.35\\
2.60205999132796	-2.35\\
2.69897000433602	-2.35\\
2.77815125038364	-2.35\\
2.84509804001426	-2.35\\
2.90308998699194	-2.35\\
2.95424250943932	-2.35\\
3	-2.35\\
};
\node[right, align=left, inner sep=0]
at (axis cs:1.35,-2.125) {\large 0.45};
\node[right, align=left, inner sep=0]
at (axis cs:2.5,-2.5) {\large 1};

\end{axis}
\end{tikzpicture}%
}}
  \caption{ Energy seminorm error assuming different AMEn
    approximation accuracies for the single edge notch tensile test
    $(a)$ and the L-shaped panel $(b)$.}
  \label{fig:energy_error_sen}
\end{figure}
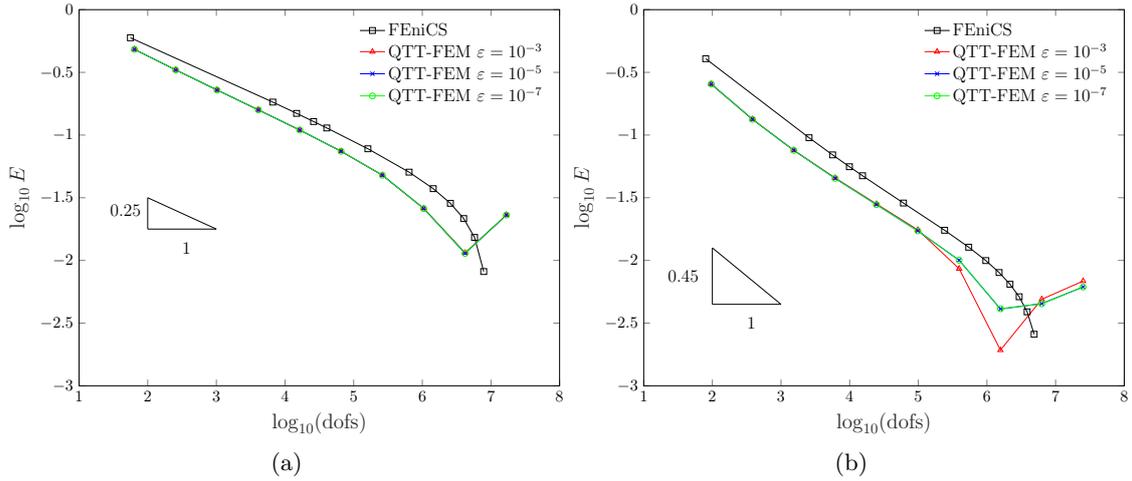
\subsection{Memory storage saving}
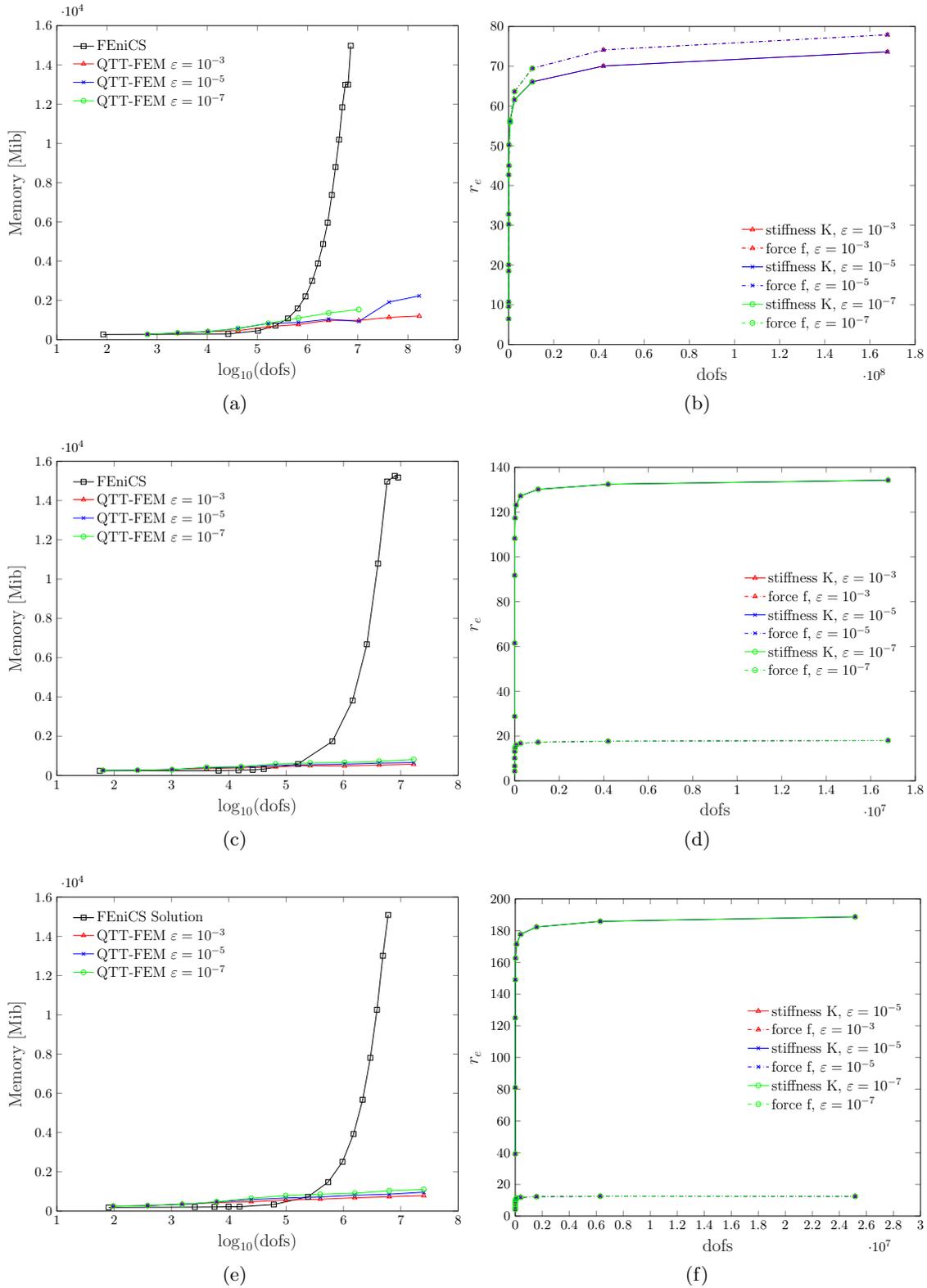
\begin{figure}[htb!]
  \centering
  \subfloat[]{\resizebox{0.45\linewidth}{!}{
%
%
\begin{tikzpicture}

\begin{axis}[%
width=4.521in,
height=3.566in,
at={(0.758in,0.481in)},
scale only axis,
xmin=1,
xmax=9,
xlabel style={font=\color{white!15!black}},
xlabel={\Large $\log_{10}$(dofs)},
ymin=0,
ymax=16000,
ylabel style={font=\color{white!15!black}},
ylabel={\Large Memory [Mib]},
axis background/.style={fill=white},
legend style={at={(0.03,0.97)}, anchor=north west, legend cell align=left, align=left, fill=none, draw=none}
]
\addplot [color=black, mark=square, mark options={solid, black}]
  table[row sep=crcr]{%
1.92427928606188	268.8359375\\
4.41584106950205	290.1171875\\
5.00903424924124	453.3710938\\
5.35822428018804	715.4804688\\
5.60659845808283	1080.449219\\
5.79951423251205	1587.589844\\
5.957272938382	2205.019531\\
6.09073477424316	2994.160156\\
6.20639460446347	3880.542969\\
6.30844744774986	4871.183594\\
6.39976057169269	5956.582031\\
6.48238072267896	7367.992188\\
6.55782011802392	8797.894531\\
6.6272277617714	10193.47266\\
6.69149718978175	11845.48047\\
6.75133700605421	12990.03906\\
6.8073186386937	13002.87109\\
6.85990961094419	14975.65625\\
};
\addlegendentry{\large  FEniCS}

\addplot [color=red, mark=triangle, mark options={solid, red}]
  table[row sep=crcr]{%
2.80617997398389	271.28125\\
3.40823996531185	322.65625\\
4.01029995663981	412.5859375\\
4.61235994796777	462.84765625\\
5.21441993929574	661.5234375\\
5.8164799306237	771.875\\
6.41853992195166	983.5625\\
7.02059991327962	980.91796875\\
7.62265990460759	1132.8046875\\
8.22471989593555	1205.140625\\
};
\addlegendentry{\large  QTT-FEM $\varepsilon=10^{-3}$}

\addplot [color=blue, mark=x, mark options={solid, blue}]
  table[row sep=crcr]{%
2.80617997398389	271.11328125\\
3.40823996531185	340.51171875\\
4.01029995663981	409.9609375\\
4.61235994796777	570.875\\
5.21441993929574	828.59375\\
5.8164799306237	871.48828125\\
6.41853992195166	1042.83984375\\
7.02059991327962	931.59375\\
7.62265990460759	1914.73828125\\
8.22471989593555	2232.3125\\
};
\addlegendentry{\large  QTT-FEM $\varepsilon=10^{-5}$}

\addplot [color=green, mark=o, mark options={solid, green}]
  table[row sep=crcr]{%
2.80617997398389	268.88671875\\
3.40823996531185	352.19140625\\
4.01029995663981	414.1171875\\
4.61235994796777	574.390625\\
5.21441993929574	833.93359375\\
5.8164799306237	1097.3828125\\
6.41853992195166	1354.890625\\
7.02059991327962	1536.77734375\\
};
\addlegendentry{\large QTT-FEM $\varepsilon=10^{-7}$}

\end{axis}
\end{tikzpicture}%
}} 
  \subfloat[]{\resizebox{0.45\linewidth}{!}{
%
%
\begin{tikzpicture}

\begin{axis}[%
width=4.521in,
height=3.566in,
at={(0.758in,0.481in)},
scale only axis,
xmin=0,
xmax=180000000,
xlabel style={font=\color{white!15!black}},
xlabel={\Large dofs},
ymin=0,
ymax=80,
ylabel style={font=\color{white!15!black}},
ylabel={\Large $r_e$},
axis background/.style={fill=white},
legend style={at={(0.97,0.03)}, anchor=south east, legend cell align=left, align=left, fill=none, draw=none}
]
\addplot [color=red, mark=triangle, mark options={solid, red}]
  table[row sep=crcr]{%
640	10.748\\
2560	20.055\\
10240	32.775\\
40960	42.696\\
163840	50.219\\
655360	56.355\\
2621440	61.554\\
10485760	66.067\\
41943040	70.05\\
167772160	73.608\\
};
\addlegendentry{\large stiffness K, $\varepsilon=10^{-3}$}

\addplot [color=red, dashdotted, mark=triangle, mark options={solid, red}]
  table[row sep=crcr]{%
640	6.472\\
2560	9.589\\
10240	18.518\\
40960	30.259\\
163840	45.026\\
655360	55.937\\
2621440	63.618\\
10485760	69.459\\
41943040	74.101\\
167772160	77.901\\
};
\addlegendentry{\large force f, $\varepsilon=10^{-3}$}

\addplot [color=blue, mark=x, mark options={solid, blue}]
  table[row sep=crcr]{%
640	10.748\\
2560	20.055\\
10240	32.775\\
40960	42.696\\
163840	50.219\\
655360	56.355\\
2621440	61.554\\
10485760	66.067\\
41943040	70.05\\
167772160	73.608\\
};
\addlegendentry{\large stiffness K,  $\varepsilon=10^{-5}$}

\addplot [color=blue, dashdotted, mark=x, mark options={solid, blue}]
  table[row sep=crcr]{%
640	6.472\\
2560	9.589\\
10240	18.518\\
40960	30.259\\
163840	45.026\\
655360	55.937\\
2621440	63.618\\
10485760	69.459\\
41943040	74.101\\
167772160	77.901\\
};
\addlegendentry{\large force f, $\varepsilon=10^{-5}$}

\addplot [color=green, mark=o, mark options={solid, green}]
  table[row sep=crcr]{%
640	10.748\\
2560	20.055\\
10240	32.775\\
40960	42.696\\
163840	50.219\\
655360	56.355\\
2621440	61.554\\
10485760	66.067\\
};
\addlegendentry{\large stiffness K, $\varepsilon=10^{-7}$}

\addplot [color=green, dashdotted, mark=o, mark options={solid, green}]
  table[row sep=crcr]{%
640	6.472\\
2560	9.589\\
10240	18.518\\
40960	30.259\\
163840	45.026\\
655360	55.937\\
2621440	63.618\\
10485760	69.459\\
};
\addlegendentry{\large force f, $\varepsilon=10^{-7}$}

\end{axis}
\end{tikzpicture}%
}}\\
  \subfloat[]{\resizebox{.45\linewidth}{!}{
%
%
\begin{tikzpicture}

\begin{axis}[%
width=4.521in,
height=3.566in,
at={(0.758in,0.481in)},
scale only axis,
xmin=1,
xmax=8,
xlabel style={font=\color{white!15!black}},
xlabel={\Large $\log_{10}$(dofs)},
ymin=0,
ymax=16000,
ylabel style={font=\color{white!15!black}},
ylabel={\Large Memory [Mib]},
axis background/.style={fill=white},
legend style={at={(0.03,0.97)}, anchor=north west, legend cell align=left, align=left, fill=none, draw=none}
]
\addplot [color=black, mark=square, mark options={solid, black}]
  table[row sep=crcr]{%
1.7481880270062	231.453125\\
3.82229887126237	245.109375\\
4.1691452009912	261.9960938\\
4.41634089057448	284.0351563\\
4.60854742686711	327.1992188\\
5.20737042570465	581.953125\\
5.8078068843684	1736.800781\\
6.15944747506963	3818.625\\
6.40905384366752	6678.433594\\
6.60271116177951	10791.15625\\
6.76096516312391	14974.72656\\
6.89478123977332	15255.4375\\
6.9546766833321	15179.37109\\
};
\addlegendentry{\large FEniCS}

\addplot [color=red, mark=triangle, mark options={solid, red}]
  table[row sep=crcr]{%
1.80617997398389	252.96484375\\
2.40823996531185	265.71484375\\
3.01029995663981	301.78515625\\
3.61235994796777	353.8671875\\
4.21441993929574	384.6953125\\
4.8164799306237	437.78515625\\
5.41853992195166	496.44140625\\
6.02059991327962	486.015625\\
6.62265990460759	532.109375\\
7.22471989593555	570.67578125\\
};
\addlegendentry{\large QTT-FEM $\varepsilon=10^{-3}$}

\addplot [color=blue, mark=x, mark options={solid, blue}]
  table[row sep=crcr]{%
1.80617997398389	253.98046875\\
2.40823996531185	265.82421875\\
3.01029995663981	302.0703125\\
3.61235994796777	402.6953125\\
4.21441993929574	429.35546875\\
4.8164799306237	513.8203125\\
5.41853992195166	548.5390625\\
6.02059991327962	573.828125\\
6.62265990460759	627.32421875\\
7.22471989593555	655.6484375\\
};
\addlegendentry{\large QTT-FEM $\varepsilon=10^{-5}$}

\addplot [color=green, mark=o, mark options={solid, green}]
  table[row sep=crcr]{%
1.80617997398389	254.2578125\\
2.40823996531185	265.8359375\\
3.01029995663981	292.265625\\
3.61235994796777	419.51171875\\
4.21441993929574	448.1796875\\
4.8164799306237	590.45703125\\
5.41853992195166	635.02734375\\
6.02059991327962	662.90234375\\
6.62265990460759	721.0859375\\
7.22471989593555	813.5078125\\
};
\addlegendentry{\large QTT-FEM $\varepsilon=10^{-7}$}

\end{axis}
\end{tikzpicture}%
}} 
  \subfloat[]{\resizebox{.45\linewidth}{!}{
%
%
\begin{tikzpicture}

\begin{axis}[%
width=4.521in,
height=3.566in,
at={(0.758in,0.481in)},
scale only axis,
xmin=0,
xmax=18000000,
xlabel style={font=\color{white!15!black}},
xlabel={\Large dofs},
ymin=0,
ymax=140,
ylabel style={font=\color{white!15!black}},
ylabel={\Large $r_e$},
axis background/.style={fill=white},
legend style={at={(0.97,0.5)}, anchor=east, legend cell align=left, align=left, fill=none, draw=none}
]
\addplot [color=red, mark=triangle, mark options={solid, red}]
  table[row sep=crcr]{%
64	28.7754163078853\\
256	61.4786383931076\\
1024	91.7493463687287\\
4096	108.288980238555\\
16384	117.373059029512\\
65536	123.170799461997\\
262144	127.205504716664\\
1048576	130.179491768115\\
4194304	132.464138880455\\
16777216	134.274984056628\\
};
\addlegendentry{\large stiffness K, $\varepsilon=10^{-3}$}

\addplot [color=red, dashdotted, mark=triangle, mark options={solid, red}]
  table[row sep=crcr]{%
64	4.35626742811115\\
256	6.66347614342147\\
1024	10.1802281448299\\
4096	13.1593730050937\\
16384	14.8793994310221\\
65536	15.9665524018049\\
262144	16.7209195150958\\
1048576	17.2765911441694\\
4194304	17.7035372443311\\
16777216	18.0421116716065\\
};
\addlegendentry{\large force f, $\varepsilon=10^{-3}$}

\addplot [color=blue, mark=x, mark options={solid, blue}]
  table[row sep=crcr]{%
64	28.7754163078853\\
256	61.4786383931076\\
1024	91.7493463687287\\
4096	108.288980238555\\
16384	117.373059029512\\
65536	123.170799461997\\
262144	127.205504716664\\
1048576	130.179491768115\\
4194304	132.464138880455\\
16777216	134.274984056628\\
};
\addlegendentry{\large stiffness K, $\varepsilon=10^{-5}$}

\addplot [color=blue, dashdotted, mark=x, mark options={solid, blue}]
  table[row sep=crcr]{%
64	4.35626742811115\\
256	6.66347614342147\\
1024	10.1802281448299\\
4096	13.1593730050937\\
16384	14.8793994310221\\
65536	15.9665524018049\\
262144	16.7209195150958\\
1048576	17.2765911441694\\
4194304	17.7035372443311\\
16777216	18.0421116716065\\
};
\addlegendentry{\large force f, $\varepsilon=10^{-5}$}

\addplot [color=green, mark=o, mark options={solid, green}]
  table[row sep=crcr]{%
64	28.7754163078853\\
256	61.4786383931076\\
1024	91.7493463687287\\
4096	108.288980238555\\
16384	117.373059029512\\
65536	123.170799461997\\
262144	127.205504716664\\
1048576	130.179491768115\\
4194304	132.464138880455\\
16777216	134.274984056628\\
};
\addlegendentry{\large stiffness K, $\varepsilon=10^{-7}$}

\addplot [color=green, dashdotted, mark=o, mark options={solid, green}]
  table[row sep=crcr]{%
64	4.35626742811115\\
256	6.66347614342147\\
1024	10.1802281448299\\
4096	13.1593730050937\\
16384	14.8793994310221\\
65536	15.9665524018049\\
262144	16.7209195150958\\
1048576	17.2765911441694\\
4194304	17.7035372443311\\
16777216	18.0421116716065\\
};
\addlegendentry{\large force f, $\varepsilon=10^{-7}$}

\end{axis}
\end{tikzpicture}%
}}\\
   \subfloat[]{\resizebox{.45\linewidth}{!}{
%
%
\begin{tikzpicture}

\begin{axis}[%
width=4.521in,
height=3.566in,
at={(0.758in,0.481in)},
scale only axis,
xmin=1,
xmax=8,
xlabel style={font=\color{white!15!black}},
xlabel={\Large $\log_{10}$(dofs)},
ymin=0,
ymax=16000,
ylabel style={font=\color{white!15!black}},
ylabel={\Large Memory [Mib]},
axis background/.style={fill=white},
legend style={at={(0.03,0.97)}, anchor=north west, legend cell align=left, align=left, fill=none, draw=none}
]
\addplot [color=black, mark=square, mark options={solid, black}]
  table[row sep=crcr]{%
1.90308998699194	186.39453125\\
3.40857912540867	203.30859375\\
3.75143308181935	210.03125\\
3.99659922270017	220.1289063\\
4.18757711905509	220.98046875\\
4.78391786504047	333.51953125\\
5.38310052508721	729.89453125\\
5.73432128223689	1481.01953125\\
5.98371737568707	2513.265625\\
6.1772484137745	3933.19140625\\
6.33541818048625	5673.75\\
6.46917406603659	7810.921875\\
6.58505467180525	10256.83984375\\
6.6872793708305	13010.66015625\\
6.77873006856815	15089.046875\\
};
\addlegendentry{\large FEniCS Solution}

\addplot [color=red, mark=triangle, mark options={solid, red}]
  table[row sep=crcr]{%
1.98227123303957	255.17578125\\
2.58433122436753	277.46484375\\
3.18639121569549	355.67578125\\
3.78845120702346	430.359375\\
4.39051119835142	477.515625\\
4.99257118967938	562.75\\
5.59463118100734	606.80078125\\
6.19669117233531	677.4765625\\
6.79875116366327	730.6796875\\
7.40081115499123	790.40625\\
};
\addlegendentry{\large QTT-FEM $\varepsilon=10^{-3}$}

\addplot [color=blue, mark=x, mark options={solid, blue}]
  table[row sep=crcr]{%
1.98227123303957	255.16015625\\
2.58433122436753	276.9609375\\
3.18639121569549	341.08203125\\
3.78845120702346	465.55859375\\
4.39051119835142	581.56640625\\
4.99257118967938	668.0859375\\
5.59463118100734	717.91015625\\
6.19669117233531	805.4140625\\
6.79875116366327	861.1328125\\
7.40081115499123	963.3671875\\
};
\addlegendentry{\large QTT-FEM $\varepsilon=10^{-5}$}

\addplot [color=green, mark=o, mark options={solid, green}]
  table[row sep=crcr]{%
1.98227123303957	256.03125\\
2.58433122436753	277.0703125\\
3.18639121569549	348.73828125\\
3.78845120702346	474.9765625\\
4.39051119835142	650.98046875\\
4.99257118967938	786.44140625\\
5.59463118100734	856.08984375\\
6.19669117233531	911.78125\\
6.79875116366327	1033.875\\
7.40081115499123	1104.859375\\
};
\addlegendentry{\large QTT-FEM $\varepsilon=10^{-7}$}

\end{axis}
\end{tikzpicture}%
}} 
  \subfloat[]{\resizebox{.45\linewidth}{!}{
%
%
\begin{tikzpicture}

\begin{axis}[%
width=4.521in,
height=3.566in,
at={(0.758in,0.481in)},
scale only axis,
xmin=0,
xmax=30000000,
xlabel style={font=\color{white!15!black}},
xlabel={\Large dofs},
ymin=0,
ymax=200,
ylabel style={font=\color{white!15!black}},
ylabel={\Large $r_e$},
axis background/.style={fill=white},
legend style={at={(0.97,0.5)}, anchor=east, legend cell align=left, align=left, fill=none, draw=none}
]
\addplot [color=red, mark=triangle, mark options={solid, red}]
  table[row sep=crcr]{%
96	39.2\\
384	81.077267370143\\
1536	125.024314386979\\
6144	149.11729565749\\
24576	162.691207329017\\
98304	171.498994030629\\
393216	177.700230339619\\
1572864	182.310911577068\\
6291456	185.876662788985\\
25165824	188.718030160536\\
};
\addlegendentry{\large stiffness K, $\varepsilon=10^{-5}$}

\addplot [color=red, dashdotted, mark=triangle, mark options={solid, red}]
  table[row sep=crcr]{%
96	4.28138329195272\\
384	5.83850070845756\\
1536	7.74397045983094\\
6144	9.3191514502424\\
24576	10.4984661470097\\
98304	11.271678583624\\
393216	11.8145789196152\\
1572864	12.2178939456219\\
6291456	12.5297812538225\\
25165824	12.4334500996882\\
};
\addlegendentry{\large force f, $\varepsilon=10^{-3}$}

\addplot [color=blue, mark=x, mark options={solid, blue}]
  table[row sep=crcr]{%
96	39.2\\
384	81.077267370143\\
1536	125.024314386979\\
6144	149.11729565749\\
24576	162.691207329017\\
98304	171.498994030629\\
393216	177.700230339619\\
1572864	182.310911577068\\
6291456	185.876662788985\\
25165824	188.718030160536\\
};
\addlegendentry{\large stiffness K, $\varepsilon=10^{-5}$}

\addplot [color=blue, dashdotted, mark=x, mark options={solid, blue}]
  table[row sep=crcr]{%
96	4.2813832919527\\
384	5.8385007084576\\
1536	7.7439704598309\\
6144	9.3191514502424\\
24576	10.4984661470097\\
98304	11.271678583624\\
393216	11.8145789196152\\
1572864	12.2178939456219\\
6291456	12.5297812538225\\
25165824	12.4334500996882\\
};
\addlegendentry{\large force f, $\varepsilon=10^{-5}$}

\addplot [color=green, mark=o, mark options={solid, green}]
  table[row sep=crcr]{%
96	39.2\\
384	81.0772673701431\\
1536	125.024314386979\\
6144	149.11729565749\\
24576	162.691207329017\\
98304	171.498994030629\\
393216	177.700230339619\\
1572864	182.310911577068\\
6291456	185.876662788985\\
25165824	188.718030160536\\
};
\addlegendentry{\large stiffness K, $\varepsilon=10^{-7}$}

\addplot [color=green, dashdotted, mark=o, mark options={solid, green}]
  table[row sep=crcr]{%
96	4.2813832919527\\
384	5.8385007084576\\
1536	7.7439704598309\\
6144	9.3191514502424\\
24576	10.4984661470097\\
98304	11.271678583624\\
393216	11.8145789196152\\
1572864	12.2178939456219\\
6291456	12.5297812538225\\
25165824	12.4334500996882\\
};
\addlegendentry{\large force f, $\varepsilon=10^{-7}$}

\end{axis}
\end{tikzpicture}%
}}
  \caption{Cantilever beam test case: $(a)$ memory storage versus the
    number of degrees of freedom for different approximation accuracy
    $\varepsilon$ and $(b)$ effective rank $r_e$ versus the
    number of degrees of freedom.
    Single edge notch tensile test: $(c)$ memory storage versus
    degrees of freedom for different values of $\varepsilon$, and $(d)$ effective rank $r_e$ versus the number of degrees of
    freedom.
    L-shaped panel: $(e)$ memory storage versus the number of degrees
    of freedom for different approximation accuracy $\varepsilon$, and
    $(f)$ effective rank $r_e$ versus the number of degrees of
    freedom.}
  \label{fig:memory}
\end{figure}

Figures \FIGREF{fig:memory}{a}, \FIGREF{fig:memory}{c}, and
\FIGREF{fig:memory}{e} highlight the drastic memory cut of the QTT
memory storage compared to the memory storage required by \FEniCS{} for
increasing degrees of freedom.
\FEniCS{} results exhibit indeed a vertiginous peak of memory
consumption.
Remarkably, this memory-peak-cut ability is quite stable when assessed
for different approximation accuracies $\varepsilon$ of the AMEn solver.

The effective rank $r_e$ (erank) is also plotted as a function of the dofs in Figures \FIGREF{fig:memory}{b}, \FIGREF{fig:memory}{d} and
\FIGREF{fig:memory}{f}.
The exhibited asymptotic trend of $r_e$ is a consequence of adopting the Z-order-based Kronecker product between matrices \cite{markeeva2020}.
Using the Z-ordering system instead of the canonical order is clearly a major improvement since the latter cannot achieve such an advantageous bound of the effective rank.

\subsection{Exponential convergence}

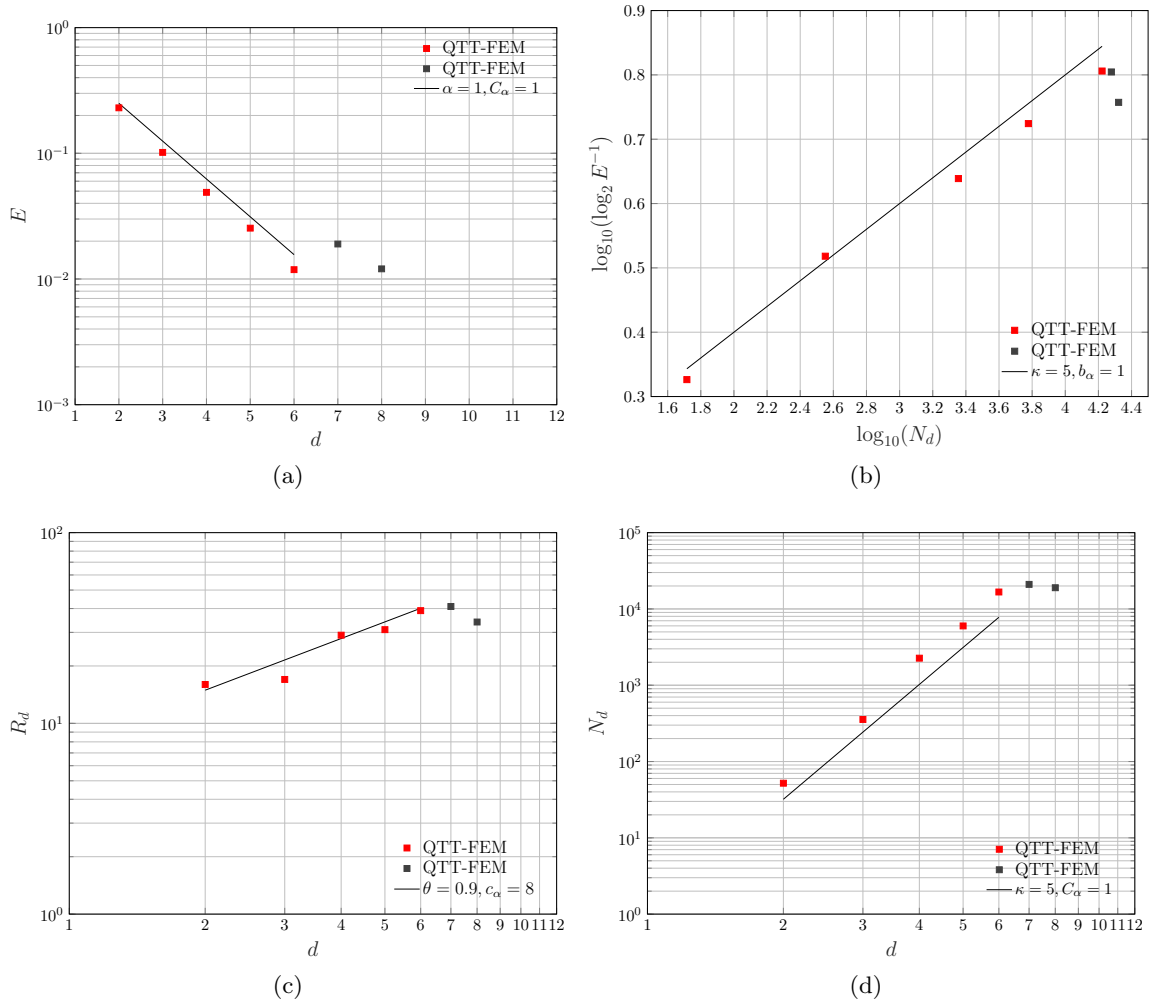
\begin{figure}[htb!]
  \centering
  \subfloat[]{\resizebox{.45\linewidth}{!}{
%
%
\definecolor{mycolor1}{rgb}{0.00000,0.44700,0.74100}%
\definecolor{mycolor2}{rgb}{0.85000,0.32500,0.09800}%
\begin{tikzpicture}

\begin{axis}[%
width=4.521in,
height=3.566in,
at={(0.758in,0.481in)},
scale only axis,
xmin=1,
xmax=12,
xlabel style={font=\color{white!15!black}},
xlabel={\Large $d$},
ymode=log,
ymin=0.001,
ymax=1,
yminorticks=true,
ylabel style={font=\color{white!15!black}},
ylabel={\Large $E$},
axis background/.style={fill=white},
xmajorgrids,
ymajorgrids,
yminorgrids,
legend style={legend cell align=left, align=left, fill=none, draw=none}
]
\addplot [color=mycolor1, only marks, mark=square*, mark options={solid, fill=red, draw=red}]
  table[row sep=crcr]{%
2	0.230071426451362\\
3	0.101722999710205\\
4	0.0489272113220465\\
5	0.0253873414447273\\
6	0.0118641247925333\\
};
\addlegendentry{\large QTT-FEM}

\addplot [color=mycolor2, only marks, mark=square*, mark options={solid, fill=darkgray, draw=darkgray}]
  table[row sep=crcr]{%
7	0.0189843316064749\\
8	0.0120416526497181\\
};
\addlegendentry{\large QTT-FEM}

\addplot [color=black]
  table[row sep=crcr]{%
2	0.25\\
3	0.125\\
4	0.0625\\
5	0.03125\\
6	0.015625\\
};
\addlegendentry{$\alpha=1, C_{\alpha}=1$}

\end{axis}
\end{tikzpicture}%
}}\hspace{0.5mm}
  \subfloat[]{\resizebox{.45\linewidth}{!}{
%
%
\definecolor{mycolor1}{rgb}{0.00000,0.44700,0.74100}%
\definecolor{mycolor2}{rgb}{0.85000,0.32500,0.09800}%
\begin{tikzpicture}

\begin{axis}[%
width=4.521in,
height=3.538in,
at={(0.758in,0.509in)},
scale only axis,
xmin=1.5,
xmax=4.5,
xlabel style={font=\color{white!15!black}},
xlabel={\Large $\log_{10}(N_{d})$},
ymin=0.3,
ymax=0.9,
ylabel style={font=\color{white!15!black}},
ylabel={\Large $\log_{10}(\log_{2} E^{-1})$},
axis background/.style={fill=white},
xmajorgrids,
ymajorgrids,
legend style={at={(0.97,0.03)}, anchor=south east, legend cell align=left, align=left, fill=none, draw=none}
]
\addplot [color=mycolor1, only marks, mark=square*, mark options={solid, fill=red, draw=red}]
  table[row sep=crcr]{%
1.7160033436348	0.326304368240023\\
2.55144999797288	0.518156115927603\\
3.35487642251623	0.638810529306009\\
3.77728179167101	0.724255126668291\\
4.22209188645426	0.805993357717239\\
};
\addlegendentry{\large QTT-FEM}

\addplot [color=mycolor2, only marks, mark=square*, mark options={solid, fill=darkgray, draw=darkgray}]
  table[row sep=crcr]{%
4.3219710555263	0.757323664430312\\
4.27811311597983	0.804536237220027\\
};
\addlegendentry{\large QTT-FEM}

\addplot [color=black]
  table[row sep=crcr]{%
1.7160033436348	0.34320066872696\\
2.55144999797288	0.510289999594575\\
3.35487642251623	0.670975284503247\\
3.77728179167101	0.755456358334203\\
4.22209188645426	0.844418377290852\\
};
\addlegendentry{$\kappa=5, b_{\alpha}=1$}

\end{axis}
\end{tikzpicture}%
}}\\
  \subfloat[]{\resizebox{.45\linewidth}{!}{
%
%
\definecolor{mycolor1}{rgb}{0.00000,0.44700,0.74100}%
\definecolor{mycolor2}{rgb}{0.85000,0.32500,0.09800}%
\begin{tikzpicture}

\begin{axis}[%
width=4.521in,
height=3.566in,
at={(0.758in,0.481in)},
scale only axis,
xmode=log,
xmin=1,
xmax=12,
xtick={1,2,3,4,5,6,7,8,9,10,11,12},
xticklabels={{1},{2},{3},{4},{5},{6},{7},{8},{9},{10},{11},{12}},
xminorticks=true,
xlabel style={font=\color{white!15!black}},
xlabel={\Large $d$},
ymode=log,
ymin=1,
ymax=100,
yminorticks=true,
ylabel style={font=\color{white!15!black}},
ylabel={\Large $R_{d}$},
axis background/.style={fill=white},
xmajorgrids,
xminorgrids,
ymajorgrids,
yminorgrids,
legend style={at={(0.97,0.03)}, anchor=south east, legend cell align=left, align=left, fill=none, draw=none}
]
\addplot [color=mycolor1, only marks, mark=square*, mark options={solid, fill=red, draw=red}]
  table[row sep=crcr]{%
2	16\\3	17\\4	29\\5	31\\6	39\\
};
\addlegendentry{\large QTT-FEM}

\addplot [color=mycolor2, only marks, mark=square*, mark options={solid, fill=darkgray, draw=darkgray}]
  table[row sep=crcr]{%
7	41\\
8	34\\
};
\addlegendentry{\large QTT-FEM}

\addplot [color=black]
  table[row sep=crcr]{%
2	14.9285278645889\\
3	21.5030030361783\\
4	27.857618025476\\
5	34.0535969008314\\
6	40.126022499741\\
};
\addlegendentry{\large $\theta=0.9, c_{\alpha}=8$}

\end{axis}
\end{tikzpicture}%
}}\hspace{0.5mm}
  \subfloat[]{\resizebox{.45\linewidth}{!}{
%
%
\definecolor{mycolor1}{rgb}{0.00000,0.44700,0.74100}%
\definecolor{mycolor2}{rgb}{0.85000,0.32500,0.09800}%
\begin{tikzpicture}

\begin{axis}[%
width=4.521in,
height=3.566in,
at={(0.758in,0.481in)},
scale only axis,
xmode=log,
xmin=1,
xmax=12,
xtick={1,2,3,4,5,6,7,8,9,10,11,12},
xticklabels={{1},{2},{3},{4},{5},{6},{7},{8},{9},{10},{11},{12}},
xminorticks=true,
xlabel style={font=\color{white!15!black}},
xlabel={\Large $d$ },
ymode=log,
ymin=1,
ymax=100000,
yminorticks=true,
ylabel style={font=\color{white!15!black}},
ylabel={\Large $N_{d}$},
axis background/.style={fill=white},
xmajorgrids,
xminorgrids,
ymajorgrids,
yminorgrids,
legend style={at={(0.97,0.03)}, anchor=south east, legend cell align=left, align=left, fill=none, draw=none}
]
\addplot [color=mycolor1, only marks, mark=square*, mark options={solid, fill=red, draw=red}]
  table[row sep=crcr]{%
2	52\\
3	356\\
4	2264\\
5	5988\\
6	16676\\
};
\addlegendentry{\large QTT-FEM}

\addplot [color=mycolor2, only marks, mark=square*, mark options={solid, fill=darkgray, draw=darkgray}]
  table[row sep=crcr]{%
7	20988\\
8	18972\\
};
\addlegendentry{\large QTT-FEM}

\addplot [color=black]
  table[row sep=crcr]{%
2	32\\
3	243\\
4	1024\\
5	3125\\
6	7776\\
};
\addlegendentry{$\kappa=5, C_{\alpha}=1$}

\end{axis}
\end{tikzpicture}%
}}
  \caption{ Cantilever beam test case: $(a)$ Energy seminorm error $E$
    (red dots) and reference exponential convergence of equation
    $C_{\alpha}2^{-\alpha d}$ (continuous line) for increasing number
    of levels $d$.
    $(b)$ Energy seminorm error $E$ (red dots) and reference line
    $\log^{k}_{2} E_{d}=-b_{\alpha}N_{d}$ for variable number of
    parameters $N_{d}$.
    $(c)$ Maximum rank of QTT approximation $R_{d}$ (red dots) and
    reference line $R_{d}=c_{\alpha}d^{\theta}$ (continuous line)
    versus $d$ levels.
    $(d)$ Number of parameters $N_{d}$ (red dots) and reference line
    $N_{d}=C_{\alpha}d^{\kappa}$ (continuous line) for increasing $d$
    levels.
    We report the points of the numerical tests with significant
    rounding error in gray, while the red points indicate the expected
    convergence rate.}
    \label{fig:kazeev_cantilever}
\end{figure}

\begin{figure}[htb!]
  \centering
  \subfloat[]{\resizebox{.45\linewidth}{!}{
%
%
\definecolor{mycolor1}{rgb}{0.00000,0.44700,0.74100}%
\definecolor{mycolor2}{rgb}{0.85000,0.32500,0.09800}%
\begin{tikzpicture}

\begin{axis}[%
width=4.521in,
height=3.566in,
at={(0.758in,0.481in)},
scale only axis,
xmin=1,
xmax=12,
xlabel style={font=\color{white!15!black}},
xlabel={\Large $d$},
ymode=log,
ymin=0.001,
ymax=1,
yminorticks=true,
ylabel style={font=\color{white!15!black}},
ylabel={\Large $E$},
axis background/.style={fill=white},
xmajorgrids,
ymajorgrids,
yminorgrids,
legend style={legend cell align=left, align=left, fill=none, draw=none}
]
\addplot [color=mycolor1, only marks, mark=square*, mark options={solid, fill=red, draw=red}]
  table[row sep=crcr]{%
2	0.482606177929383\\
3	0.330640940445522\\
4	0.228732940132761\\
5	0.158750563340338\\
6	0.109632616919392\\
7	0.0742905538705022\\
8	0.0478291886897279\\
9	0.0259831111352536\\
10	0.0115047358128803\\
};
\addlegendentry{\large QTT-FEM}

\addplot [color=mycolor2, only marks, mark=square*, mark options={solid, fill=darkgray, draw=darkgray}]
  table[row sep=crcr]{%
11	0.0229809828373837\\
};
\addlegendentry{\large QTT-FEM}

\addplot [color=black]
  table[row sep=crcr]{%
2	0.5\\
3	0.353553390593274\\
4	0.25\\
5	0.176776695296637\\
6	0.125\\
7	0.0883883476483184\\
8	0.0625\\
9	0.0441941738241592\\
10	0.03125\\
};
\addlegendentry{\large $\alpha\text{=0.5}$, $C_{\alpha}=1$}

\end{axis}
\end{tikzpicture}%
}}\hspace{0.5mm}
  \subfloat[]{\resizebox{.45\linewidth}{!}{
%
%
\definecolor{mycolor1}{rgb}{0.00000,0.44700,0.74100}%
\definecolor{mycolor2}{rgb}{0.85000,0.32500,0.09800}%
\begin{tikzpicture}

\begin{axis}[%
width=4.521in,
height=3.538in,
at={(0.758in,0.509in)},
scale only axis,
xmin=1.5,
xmax=5,
xlabel style={font=\color{white!15!black}},
xlabel={\Large $\log_{10}(N_{d})$},
ymin=0,
ymax=0.9,
ylabel style={font=\color{white!15!black}},
ylabel={\Large $\log_{10}(\log_{2} E^{-1})$},
axis background/.style={fill=white},
xmajorgrids,
ymajorgrids,
legend style={at={(0.97,0.03)}, anchor=south east, legend cell align=left, align=left, fill=none, draw=none}
]
\addplot [color=mycolor1, only marks, mark=square*, mark options={solid, fill=red, draw=red}]
  table[row sep=crcr]{%
1.7160033436348	0.021636479361388\\
2.55144999797288	0.203213185900166\\
3.21801004298436	0.328025488667548\\
3.73335778792559	0.424091740918759\\
4.03165079355126	0.503688702409733\\
4.2723986324554	0.574109712488035\\
4.45344065929356	0.642065148831623\\
4.52260070783044	0.721504111492794\\
4.60967976584537	0.808995663552898\\
};
\addlegendentry{\large QTT-FEM}

\addplot [color=mycolor2, only marks, mark=square*, mark options={solid, fill=darkgray, draw=darkgray}]
  table[row sep=crcr]{%
4.69460519893357	0.735871500775586\\
};
\addlegendentry{\large QTT-FEM}

\addplot [color=black]
  table[row sep=crcr]{%
1.7160033436348	0.0636126669925523\\
2.55144999797288	0.230701997860167\\
3.21801004298436	0.364014006862465\\
3.73335778792559	0.46708355585071\\
4.03165079355126	0.526742156975845\\
4.2723986324554	0.574891724756673\\
4.45344065929356	0.611100130124304\\
4.52260070783044	0.624932139831681\\
4.60967976584537	0.642347951434666\\
};
\addlegendentry{$\kappa=5$, $b_{\alpha}$=0.04}

\end{axis}
\end{tikzpicture}%
}}\\
  \subfloat[]{\resizebox{.45\linewidth}{!}{
%
%
\definecolor{mycolor1}{rgb}{0.00000,0.44700,0.74100}%
\definecolor{mycolor2}{rgb}{0.85000,0.32500,0.09800}%
\begin{tikzpicture}

\begin{axis}[%
width=4.521in,
height=3.566in,
at={(0.758in,0.481in)},
scale only axis,
xmode=log,
xmin=1,
xmax=12,
xtick={1,2,3,4,5,6,7,8,9,10,11,12},
xticklabels={{1},{2},{3},{4},{5},{6},{7},{8},{9},{10},{11},{12}},
xminorticks=true,
xlabel style={font=\color{white!15!black}},
xlabel={\Large $d$},
ymode=log,
ymin=1,
ymax=100,
yminorticks=true,
ylabel style={font=\color{white!15!black}},
ylabel={\Large$R_{d}$},
axis background/.style={fill=white},
xmajorgrids,
xminorgrids,
ymajorgrids,
yminorgrids,
legend style={at={(0.97,0.03)}, anchor=south east, legend cell align=left, align=left, fill=none, draw=none}
]
\addplot [color=mycolor1, only marks, mark=square*, mark options={solid, fill=red, draw=red}]
  table[row sep=crcr]{%
2	8\\
3	16\\
4	20\\
5	35\\
6	38\\
7	41\\
8	44\\
9	45\\
10	47\\
};
\addlegendentry{\large QTT-FEM}

\addplot [color=mycolor2, only marks, mark=square*, mark options={solid, fill=darkgray, draw=darkgray}]
  table[row sep=crcr]{%
11	50\\
};
\addlegendentry{\large QTT-FEM}

\addplot [color=black]
  table[row sep=crcr]{%
2	10\\
3	15\\
4	20\\
5	25\\
6	30\\
7	35\\
8	40\\
9	45\\
10	50\\
};
\addlegendentry{\large $\theta=1$, $c_{\alpha}=1$}

\end{axis}
\end{tikzpicture}%
}}\hspace{0.5mm}
  \subfloat[]{\resizebox{.45\linewidth}{!}{
%
%
\definecolor{mycolor1}{rgb}{0.00000,0.44700,0.74100}%
\definecolor{mycolor2}{rgb}{0.85000,0.32500,0.09800}%
\begin{tikzpicture}

\begin{axis}[%
width=4.521in,
height=3.566in,
at={(0.758in,0.481in)},
scale only axis,
xmode=log,
xmin=1,
xmax=12,
xtick={1,2,3,4,5,6,7,8,9,10,11,12},
xticklabels={{1},{2},{3},{4},{5},{6},{7},{8},{9},{10},{11},{12}},
xminorticks=true,
xlabel style={font=\color{white!15!black}},
xlabel={\Large $d$},
ymode=log,
ymin=1,
ymax=100000,
yminorticks=true,
ylabel style={font=\color{white!15!black}},
ylabel={\Large$N_{d}$},
axis background/.style={fill=white},
xmajorgrids,
xminorgrids,
ymajorgrids,
yminorgrids,
legend style={at={(0.97,0.03)}, anchor=south east, legend cell align=left, align=left, fill=none, draw=none}
]
\addplot [color=mycolor1, only marks, mark=square*, mark options={solid, fill=red, draw=red}]
  table[row sep=crcr]{%
2	52\\
3	356\\
4	1652\\
5	5412\\
6	10756\\
7	18724\\
8	28408\\
9	33312\\
10	40708\\
};
\addlegendentry{\large QTT-FEM}

\addplot [color=mycolor2, only marks, mark=square*, mark options={solid, fill=darkgray, draw=darkgray}]
  table[row sep=crcr]{%
11	49500\\
};
\addlegendentry{\large QTT-FEM}

\addplot [color=black]
  table[row sep=crcr]{%
2	32\\
3	243\\
4	1024\\
5	3125\\
6	7776\\
7	16807\\
8	32768\\
9	59049\\
10	100000\\
};
\addlegendentry{\large $\kappa=5$, $C_{\alpha}=1$}

\end{axis}
\end{tikzpicture}%
}}
  \caption{Single edge-notch tensile test case: $(a)$ Energy seminorm error
    $E$ (red dots) and reference exponential convergence of equation
    $C_{\alpha}2^{-\alpha d}$ (continuous line) for increasing number
    of levels $d$.
    $(b)$ Energy seminorm error $E$ (red dots) and reference line
    $\log^{k}_{2} E_{d}=-b_{\alpha}N_{d}$ for variable number of
    parameters $N_{d}$.
    $(c)$ Maximum rank of QTT approximation $R_{d}$ (red dots) and
    reference line $R_{d}=c_{\alpha}d^{\theta}$ (continuous line)
    versus $d$ levels.
    $(d)$ Number of parameters $N_{d}$ (red dots) and reference line
    $N_{d}=C_{\alpha}d^{\kappa}$ (continuous line) for increasing $d$
    levels.
    We report the points of the numerical tests with significant
    rounding error in gray, while the red points indicate the expected
    convergence rate.}
  \label{fig:kazeev_sen}
\end{figure}
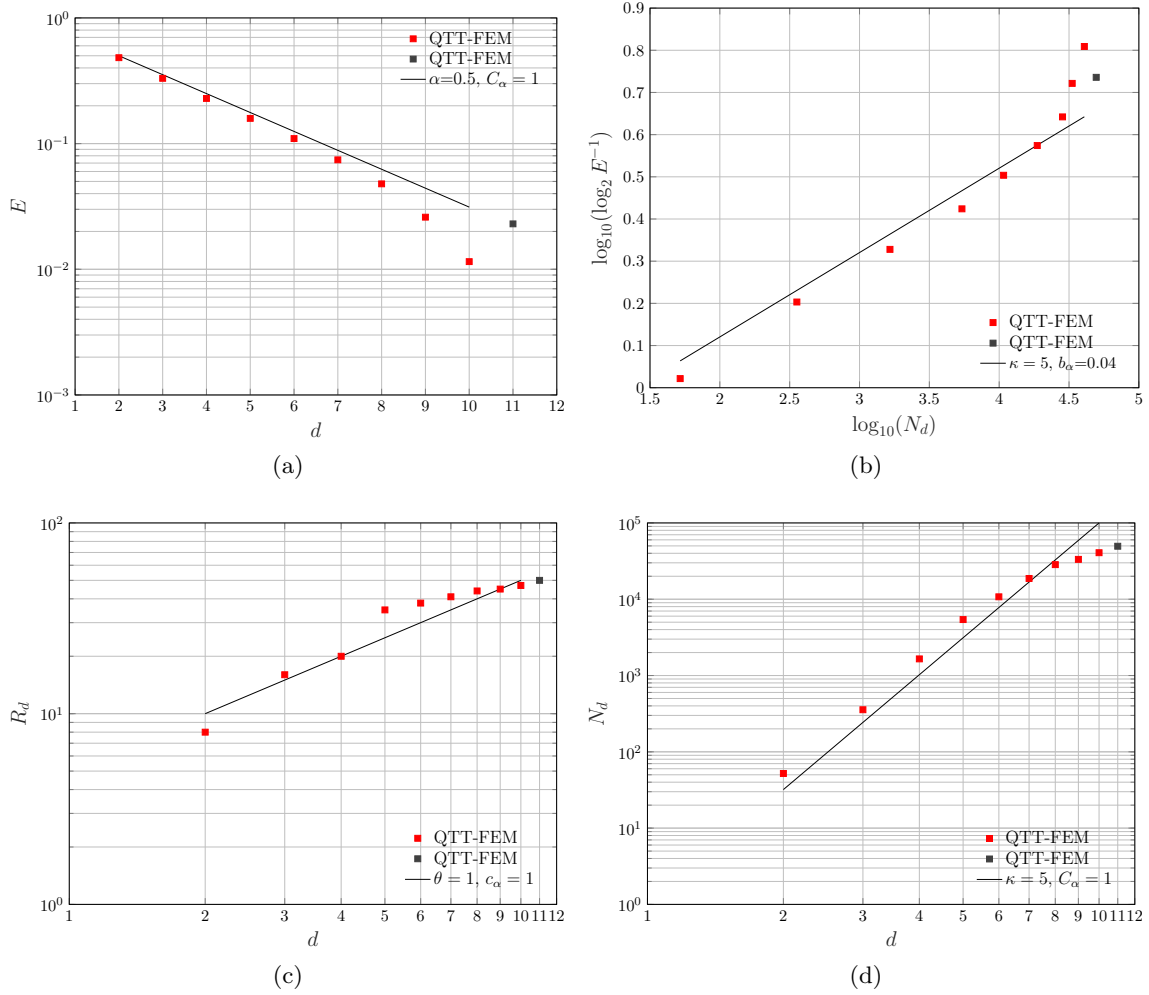

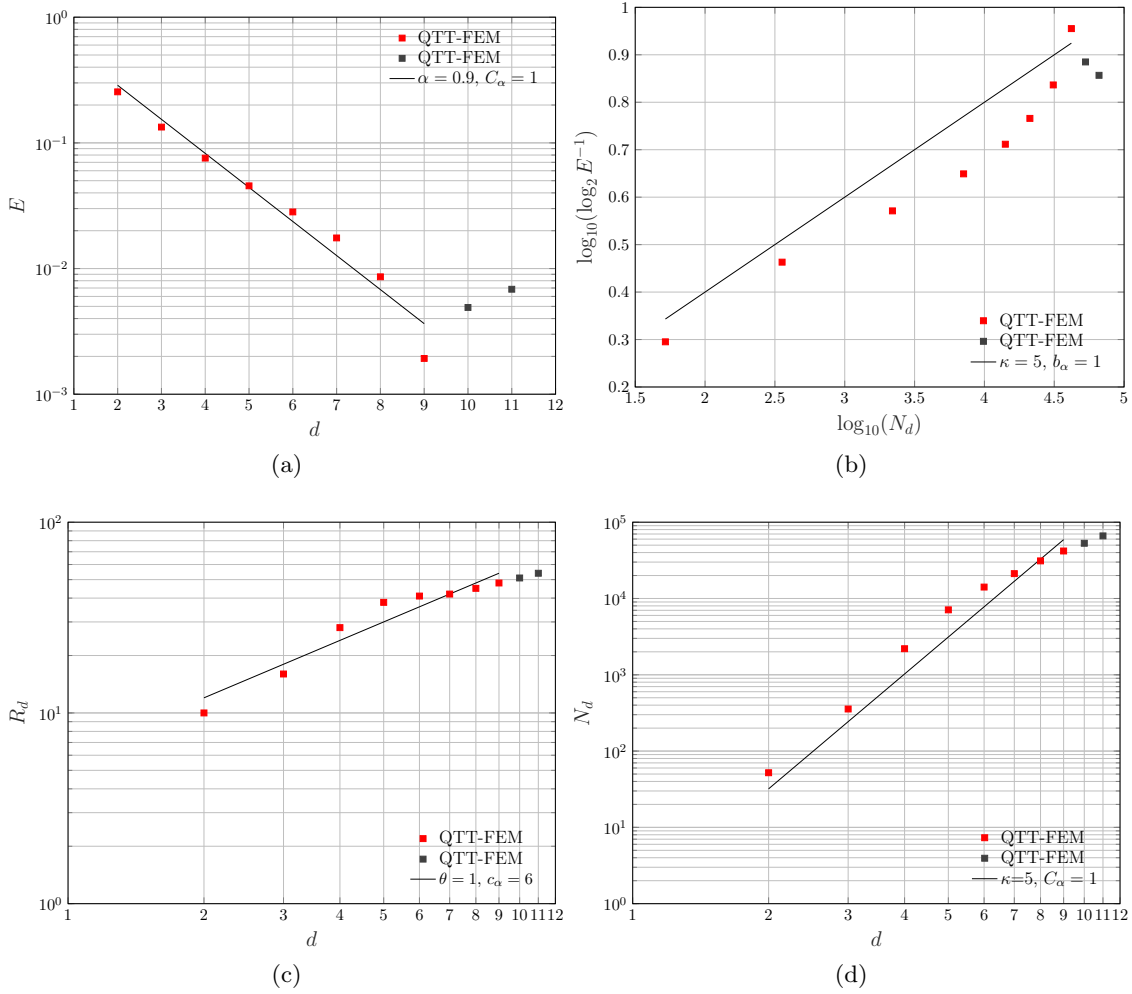
\begin{figure}[htb!]
  \centering
  \subfloat[]{\resizebox{.45\linewidth}{!}{
%
%
\definecolor{mycolor1}{rgb}{0.00000,0.44700,0.74100}%
\definecolor{mycolor2}{rgb}{0.85000,0.32500,0.09800}%
\begin{tikzpicture}

\begin{axis}[%
width=4.521in,
height=3.566in,
at={(0.758in,0.481in)},
scale only axis,
xmin=1,
xmax=12,
xlabel style={font=\color{white!15!black}},
xlabel={\Large $d$},
ymode=log,
ymin=0.001,
ymax=1,
yminorticks=true,
ylabel style={font=\color{white!15!black}},
ylabel={\Large $E$},
axis background/.style={fill=white},
xmajorgrids,
ymajorgrids,
yminorgrids,
legend style={legend cell align=left, align=left, fill=none, draw=none}
]
\addplot [color=mycolor1, only marks, mark=square*, mark options={solid, fill=red, draw=red}]
  table[row sep=crcr]{%
2	0.254675605730967\\
3	0.133551803888839\\
4	0.0756366386074908\\
5	0.0455097879808373\\
6	0.0282132521186793\\
7	0.0175254002443938\\
8	0.00859486493251946\\
9	0.00192710189246116\\
};
\addlegendentry{\large QTT-FEM }

\addplot [color=mycolor2, only marks, mark=square*, mark options={solid, fill=darkgray, draw=darkgray}]
  table[row sep=crcr]{%
10	0.0049018216921704\\
11	0.00684532023101683\\
};
\addlegendentry{\large QTT-FEM}

\addplot [color=black]
  table[row sep=crcr]{%
2	0.287174588749259\\
3	0.153893051668115\\
4	0.0824692444233059\\
5	0.0441941738241592\\
6	0.023683071351725\\
7	0.0126914436930662\\
8	0.00680117627575097\\
9	0.00364466012319065\\
};
\addlegendentry{\large $\alpha=0.9$, $C_{\alpha}=1$}

\end{axis}
\end{tikzpicture}%
}} 
  \subfloat[]{\resizebox{.45\linewidth}{!}{
%
%
\definecolor{mycolor1}{rgb}{0.00000,0.44700,0.74100}%
\definecolor{mycolor2}{rgb}{0.85000,0.32500,0.09800}%
\begin{tikzpicture}

\begin{axis}[%
width=4.521in,
height=3.538in,
at={(0.758in,0.509in)},
scale only axis,
xmin=1.5,
xmax=5,
xlabel style={font=\color{white!15!black}},
xlabel={\Large $\log_{10}(N_{d})$},
ymin=0.2,
ymax=1,
ylabel style={font=\color{white!15!black}},
ylabel={\Large $\log_{10}(\log_{2} E^{-1})$},
axis background/.style={fill=white},
xmajorgrids,
ymajorgrids,
legend style={at={(0.97,0.03)}, anchor=south east, legend cell align=left, align=left, fill=none, draw=none}
]
\addplot [color=mycolor1, only marks, mark=square*, mark options={solid, fill=red, draw=red}]
  table[row sep=crcr]{%
1.7160033436348	0.295185923027359\\
2.55144999797288	0.463075662019777\\
3.34163233577805	0.571099570398721\\
3.85101360682367	0.649108823104254\\
4.14971164968793	0.711594939220231\\
4.32559775286002	0.765996855371929\\
4.4929279301394	0.836470289858726\\
4.62254557712451	0.95517530914355\\
};
\addlegendentry{\large QTT-FEM}

\addplot [color=mycolor2, only marks, mark=square*, mark options={solid, fill=darkgray, draw=darkgray}]
  table[row sep=crcr]{%
4.72388236808015	0.884934988339703\\
4.82051671772844	0.856769131783259\\
};
\addlegendentry{\large QTT-FEM}

\addplot [color=black]
  table[row sep=crcr]{%
1.7160033436348	0.34320066872696\\
2.55144999797288	0.510289999594575\\
3.34163233577805	0.668326467155611\\
3.85101360682367	0.770202721364734\\
4.14971164968793	0.829942329937586\\
4.32559775286002	0.865119550572004\\
4.4929279301394	0.898585586027881\\
4.62254557712451	0.924509115424903\\
};
\addlegendentry{\large $\kappa=5$, $b_{\alpha}=1$}

\end{axis}
\end{tikzpicture}%
}}\\
  \subfloat[]{\resizebox{.45\linewidth}{!}{
%
%
\definecolor{mycolor1}{rgb}{0.00000,0.44700,0.74100}%
\definecolor{mycolor2}{rgb}{0.85000,0.32500,0.09800}%
\begin{tikzpicture}

\begin{axis}[%
width=4.521in,
height=3.566in,
at={(0.758in,0.481in)},
scale only axis,
xmode=log,
xmin=1,
xmax=12,
xtick={1,2,3,4,5,6,7,8,9,10,11,12},
xticklabels={{1},{2},{3},{4},{5},{6},{7},{8},{9},{10},{11},{12}},
xminorticks=true,
xlabel style={font=\color{white!15!black}},
xlabel={\Large $d$},
ymode=log,
ymin=1,
ymax=100,
yminorticks=true,
ylabel style={font=\color{white!15!black}},
ylabel={\Large $R_{d}$},
axis background/.style={fill=white},
xmajorgrids,
xminorgrids,
ymajorgrids,
yminorgrids,
legend style={at={(0.97,0.03)}, anchor=south east, legend cell align=left, align=left, fill=none, draw=none}
]
\addplot [color=mycolor1, only marks, mark=square*, mark options={solid, fill=red, draw=red}]
  table[row sep=crcr]{%
2	10\\
3	16\\
4	28\\
5	38\\
6	41\\
7	42\\
8	45\\
9	48\\
};
\addlegendentry{\large QTT-FEM}

\addplot [color=mycolor2, only marks, mark=square*, mark options={solid, fill=darkgray, draw=darkgray}]
  table[row sep=crcr]{%
10	51\\
11	54\\
};
\addlegendentry{\large QTT-FEM}

\addplot [color=black]
  table[row sep=crcr]{%
2	12\\
3	18\\
4	24\\
5	30\\
6	36\\
7	42\\
8	48\\
9	54\\
};
\addlegendentry{$\theta=1$, $c_{\alpha}=6$}

\end{axis}
\end{tikzpicture}%
}}
  \subfloat[]{\resizebox{.45\linewidth}{!}{
%
%
\definecolor{mycolor1}{rgb}{0.00000,0.44700,0.74100}%
\definecolor{mycolor2}{rgb}{0.85000,0.32500,0.09800}%
\begin{tikzpicture}

\begin{axis}[%
width=4.521in,
height=3.566in,
at={(0.758in,0.481in)},
scale only axis,
xmode=log,
xmin=1,
xmax=12,
xtick={1,2,3,4,5,6,7,8,9,10,11,12},
xticklabels={{1},{2},{3},{4},{5},{6},{7},{8},{9},{10},{11},{12}},
xminorticks=true,
xlabel style={font=\color{white!15!black}},
xlabel={\Large $d$},
ymode=log,
ymin=1,
ymax=100000,
yminorticks=true,
ylabel style={font=\color{white!15!black}},
ylabel={\Large $N_{d}$},
axis background/.style={fill=white},
xmajorgrids,
xminorgrids,
ymajorgrids,
yminorgrids,
legend style={at={(0.97,0.03)}, anchor=south east, legend cell align=left, align=left, fill=none, draw=none}
]
\addplot [color=mycolor1, only marks, mark=square*, mark options={solid, fill=red, draw=red}]
  table[row sep=crcr]{%
2	52\\
3	356\\
4	2196\\
5	7096\\
6	14116\\
7	21164\\
8	31112\\
9	41932\\
};
\addlegendentry{\large QTT-FEM}

\addplot [color=mycolor2, only marks, mark=square*, mark options={solid, fill=darkgray, draw=darkgray}]
  table[row sep=crcr]{%
10	52952\\
11	66148\\
};
\addlegendentry{\large QTT-FEM}

\addplot [color=black]
  table[row sep=crcr]{%
2	32\\
3	243\\
4	1024\\
5	3125\\
6	7776\\
7	16807\\
8	32768\\
9	59049\\
};
\addlegendentry{\large $\kappa\text{=5}$, $C_{\alpha}=1$}

\end{axis}
\end{tikzpicture}%
}}
  \caption{ L-shaped panel test case: $(a)$ Energy seminorm error $E$
    (red dots) and reference exponential convergence of equation
    $C_{\alpha}2^{-\alpha d}$ (continuous line) for increasing number
    of levels $d$.
    $(b)$ Energy seminorm error $E$ (red dots) and reference line
    $\log^{k}_{2} E_{d}=-b_{\alpha}N_{d}$ for variable number of
    parameters $N_{d}$.
    $(c)$ Maximum rank of QTT approximation $R_{d}$ (red dots) and
    reference line $R_{d}=c_{\alpha}d^{\theta}$ (continuous line)
    versus $d$ levels.
    $(d)$ Number of parameters $N_{d}$ (red dots) and reference line
    $N_{d}=C_{\alpha}d^{\kappa}$ (continuous line) for increasing $d$
    levels.
    We report the points of the numerical tests with a significant
    rounding error in gray, while the red points indicate the expected
    convergence rate.
  }
  \label{fig:kazeev_l-shaped}
\end{figure}

We further assessed the QTT convergence properties for all the examples here reported by assuming a tolerance $\varepsilon =
10^{-3}$.
Precisely, we expect an exponential convergence in the vector variable
case similar to the one proved in Ref.~\cite{kazeev2018}, where the
nodal variables are scalar quantities associated with the finite
element approximation of a Poisson problem.
We report the results in Figures \ref{fig:kazeev_cantilever},
\ref{fig:kazeev_sen}, and \ref{fig:kazeev_l-shaped} for the cantilever
beam test case, the SEN test case, and the L-shaped plate test case,
respectively.
In all such figures, the red points identify the results exhibiting
the proper energy-seminorm convergence.
For completeness, the grey points indicate the results that mostly
appear to be affected by rounding errors.

According to~\cite{kazeev2018}, we display the convergence in the
energy-seminorm error $E$ for an increasing number of levels $d$.
In this case, the continuous lines indicate the reference exponential
convergence line defined by $E=C_{\alpha}2^{-{\alpha}d}$, where
$C_{\alpha}$ is a constant factor independent of $d$ and may change
depending on the singularity order of the exact solution, and $\alpha$
is also determined by the low-order finite element approximation.
Interestingly, our results point out that
$\alpha=\textrm{min}(1,\beta)$, where $\beta$ is twice the singularity
order of the problem solution, namely $\beta=2,\,0.5,\,0.9$ for the
cantilever beam, the SEN plate, and the L-shaped panel test cases,
respectively.

We plot the error $E$ in the energy semi-norm $E$ or a variable
number of parameters $N_{d}$.
We use the line defined by $\log^{k}_{2} E =-b_{\alpha}N_{d}$ as the
reference line to highlight the exponential convergence.
The parameters $b_{\alpha}$ and $\kappa$ are independent of $d$.
The maximum rank of QTT approximation $R_{d}$ also displays an
exponential convergence in terms of the $d$ level number.
In this case, we consider the reference algebraic growth
$R_{d}=c_{\alpha}d^{\theta}$, where, again, $\theta$ and $c_{\alpha}$
are independent of $d$.
Finally, the number of parameters $N_{d}$ versus $d$ level number is
shown together with the reference algebraic growth
$N_{d}=C_{\alpha}d^{\kappa}$, $C_{\alpha}$ and $\kappa$ being
independent of $d$.

In the case of the cantilever beam, Figure
\FIGREF{fig:kazeev_cantilever}{a} shows that the expected slope of
convergence, which is equal to $1$, is reached.
Figure \FIGREF{fig:kazeev_cantilever}{b} illustrates the exponential
convergence of the QTT-FE approximation with $\kappa=5$, which
represents the theoretical upper limit demonstrated for QTT-structured
finite element discretization in~\cite{kazeev2018}.
The QTT rank $R_{d}$ grows sublinearly with increasing $d$ as shown in
Figure \FIGREF{fig:kazeev_cantilever}c{c}.
Finally, Figure \FIGREF{fig:kazeev_cantilever}(d) displays the
relation between the number of levels $d$ and the number of parameters
involved in the representation $N_{d}$.

The optimal convergence rates in Figure
\FIGREF{fig:kazeev_cantilever}{a},
the exponential convergence behavior within the theoretical limit
in Figure
\FIGREF{fig:kazeev_cantilever}{b},
the sublinear rank growth in Figure \FIGREF{fig:kazeev_cantilever}{c},
and the relationship between levels and parameters (Figure
\FIGREF{fig:kazeev_cantilever}{d},
all indicate that the QTT-FE method is performing as expected and
exhibits the desired efficiency gain compared to traditional designs
of the finite element method like that in \FEniCS.
In fact, the sublinear growth of the QTT rank when the number of levels
$d$ increases of Figure \FIGREF{fig:kazeev_cantilever}{c} is a clear
key indicator of the CPU cost and memory savings.
It demonstrates that the QTT format effectively captures the low-rank
structure of the solution, leading to significant reductions in
storage and computational complexity compared to the usual scaling of
the finite element methods with the number of degrees of freedom.
The relationship between $d$ and $N$ in Figure
\FIGREF{fig:kazeev_cantilever}{d} further illustrates how the number
of parameters in the QTT representation scales much more favorably
that the number of degrees of freedom in
a plain finite element formulation and nimplementation.
These results validate the effectiveness of the QTT approach for this
problem.

We can draw analogous conclusions for the SEN plate test case, see
Figures \FIGREF{fig:kazeev_sen}{a}-\FIGREF{fig:kazeev_sen}{d}, and the
L-shape panel test case, see Figures
\FIGREF{fig:kazeev_l-shaped}{a}-\FIGREF{fig:kazeev_l-shaped}{d}.
In these last cases, we note that 
the order of the singularity of the ground truth solution
dictates the energy convergence rate for increasing $d$ levels.
These results could be improved through the AMR technique or by
enriching the finite element space with additional terms that may
better represent the singular behavior of the solution.
However, an investigation of these methodologies is beyond he goals and
scope of the present work and will be
the subject of future research work.

\section{Conclusions} \label{sec6:conclusions}


We have extended the design of the Quantum Tensor Train finite element
solver proposed by Markeeva et al. \cite{markeeva2020} for scalar
problems to vector problems such as the linear elasticity.
Combined with Z-ordering and subdomain concatenation, our approach
achieves significant memory savings and remarkable rank reduction
compared to traditional Finite Element solvers, such as \FEniCS, all
while ensuring exponential convergence versus the number of degrees of
freedom.
The trade-off lies in a fundamental shift in the implementation
paradigm for essential finite element operations, including mesh
discretization, ordering of nodes and degrees of freedom, assembly of
stiffness matrices and internal nodal forces, and algebraic
matrix-vector computations.
In conclusion, our work confirms that using the QTT format
substantially reduces the memory usage and provides a notable
enhancement in computational speed compared to traditional
sparse-matrix finite element solvers.
Additionally, the rank growth is effectively constrained through the
use of Z-order operations.

\section*{Acknowledgments}
The Laboratory Directed Research and Development (LDRD) program
financially supported the work of G. Manzini. Los Alamos National
Laboratory is operated by Triad National Security, LLC, for the
National Nuclear Security Administration of the U.S. Department of
Energy (Contract No.\ 89233218CNA000001).
G. Manzini is a member of the Gruppo Nazionale Calcolo
Scientifico-Istituto Nazionale di Alta Matematica (GNCS-INdAM).
E. Benvenuti, M. Nale, and S. Pizzolato gratefully acknowledge the
financial support of the Italy-Croatia Interreg Project STRENGTH (ID:
ITHR0200318).


\newpage

\end{document}